\documentclass[12pt,reqno]{article}
\usepackage{amsfonts}
\usepackage[utf8x]{inputenc}
\usepackage[english]{babel}
\usepackage{amsmath,amssymb,amsfonts,wasysym}   
\usepackage{graphics}                           
\usepackage{amsmath}  
\usepackage{mathtools} 
\usepackage{hyperref}                           
\usepackage{comment}
\newtheorem{theorem}{Theorem}[section]
\newcommand*{\ttt}{\mathop{}\!\mathrm{t}}
\newcommand*{\uu}{\mathop{}\!\mathrm{u}}
\newtheorem{proposition}[theorem]{Proposition}
\newtheorem{corollary}[theorem]{Corollary}
\newtheorem{remark}[theorem]{Remark}
\newtheorem{lemma}[theorem]{Lemma}
\newtheorem{observation}[theorem]{Observation}
\newtheorem{example}[theorem]{Example}
\newtheorem{definition}[theorem]{Definition}
\newcommand*{\eff}{\mathop{}\!\mathrm{eff}}
\newcommand*{\dimensione}{\mathop{}\!\mathrm{dim}}

\newcommand*{\hh}{\mathop{}\!\mathrm{h}}
\newcommand*{\distanza}{\mathop{}\!\mathrm{dist}}
\newcommand*{\vv}{\mathop{}\!\mathrm{v}}
\newcommand*{\ww}{\mathop{}\!\mathrm{w}}
\newcommand*{\supp}{\mathop{}\!\mathrm{supp}}
\newcommand*{\valutation}{\mathop{}\!\mathrm{val}}

\newcommand*{\spaned}{\mathop{}\!\mathrm{span}}
\newcommand*{\kernel}{\mathop{}\!\mathrm{Ker}}
\newcommand*{\lie}{\mathop{}\!\mathrm{Lie}}
\newcommand*{\rk}{\mathop{}\!\mathrm{Rank}}

\newcommand*{\tr}{\mathop{}\!\mathrm{Tr}}
\newcommand*{\Proof}{\mathop{}\!\mathit{Proof}}

\usepackage{thmtools}
\usepackage{thm-restate}

\usepackage{hyperref}

\usepackage{cleveref}

\begin{document}
\title{Scaling asymptotics of Szeg\"{o} kernels under commuting Hamiltonian actions}
\author{Simone Camosso$^*$}
\date{}
{\renewcommand{\thefootnote}{\fnsymbol{footnote}}
\setcounter{footnote}{1}
\footnotetext{\textbf{Address}: Dipartimento di Matematica e Applicazioni, Università degli Studi di Milano Bicocca, Via R. Cozzi 53, 20125 Milano, Italy; \textbf{e-mail}: s.camosso@campus.unimib.it}
\setcounter{footnote}{0}
}

\maketitle

\begin{center}
\textbf{Abstract}
\end{center}
\noindent
Let $M$ be a connected $d_{M}$-dimensional complex projective manifold, and let $A$ be a holomorphic positive Hermitian line bundle on $M$, 
with normalized curvature $\omega$. Let $G$ be a compact and connected Lie group of dimension $d_{G}$, and let $T$ 
be a compact torus of dimension $d_{T}$. Suppose that both $G$ and $T$ act on $M$ in a holomorphic and Hamiltonian manner, 
that the actions commute, and linearize to $A$. If $X$ is the principal $S^1$-bundle associated to $A$, then this set-up determines
commuting unitary representations of $G$ and $T$ on the Hardy space $H(X)$ of $X$, which may then be decomposed over the irreducible 
representations of the two groups. If the moment map for the $T$-action is nowhere zero, all isotypical components for the torus  are finite-dimensional, 
and thus provide a collection of finite-dimensional $G$-modules. Given a non-zero integral weight $\nu_{T}$ for $T$, we consider the isotypical components associated to the multiples $k\nu_{T}$, $k\rightarrow +\infty$, and focus on how their structure as $G$-modules is reflected by certain local scaling asymptotics  on $X$ (and $M$). More precisely, given a fixed irreducible character $\nu_{G}$ of $G$, we study the local scaling asymptotics of the equivariant Szeg\"{o} projectors associated to $\nu_{G}$ and $k\nu_{T}$, for $k\rightarrow +\infty$, investigating their asymptotic concentration along certain loci defined by the moment maps.

\vspace{1cm}

\smallskip
\noindent \textbf{Keywords.} Scaling asymptotics, Tian-Yau-Zelditch expansion, Hardy space, Szeg\"{o} kernel, Hamiltonian action.
\newline
\smallskip
\noindent \textbf{AMS Subject Classification.} 53D05, 53D20, 53D50, 30H10, 32T15

\section{Introduction}

\subsection{The main Theorem}

Let $M$ be a connected $d_{M}$-dimensional, compact, complex projective manifold and $(A,h)$ be an ample positive Hermitian line bundle on $M$. 
We may assume that the curvature form of the unique compatible connection $\nabla_{A}$ is $\Theta=-2i\omega$, where $\omega$ is a 
K\"{a}hler form. Let $dV_{M}$  be the volume form $\frac{\omega^{\wedge d_{M}}}{d_{M}!}$ associated with $(M,\omega)$.

We put $h=g-i\omega$ where $g$ is the induced Riemannian structure. Suppose given two connected compact Lie groups $G$ and $T$, with $T$ a torus, of dimension $d_{G}$ and $d_{T}$, respectively, and commuting holomorphic and Hamiltonian actions $\mu^{G}:G\times M\rightarrow M$ and $\mu^{T}:T\times M \rightarrow M$. 
Thus $\mu_{g}^{G}\circ\mu_{t}^{T}=\mu_{t}^{T}\circ\mu_{g}^{G}$, $\forall (g,t)\in G\times T$ with moment maps $\Phi_{G},\Phi_{T}$.
Assume that both actions unitarily linearize to $A$, that is, that they admit metric preserving lifting $\tilde{\mu}^{G},\tilde{\mu}^{T}$.
Let $\widehat{G}$ be the collection of irreducible characters of $G$ and for any 
$\nu_{G}\in\widehat{G}$ let $V_{\nu_{G}}$ be the corresponding irreducible unitary representation. 
The action of $G$ on $A$ dualizes to an action on the dual line bundle $A^{\vee}$ and the $G$-invariant Hermitian metric $h$ on $A$ naturally induces an Hermitian metric on $A^{\vee}$ also denoted by $h$.

Let $X\subseteq A^{\vee}$ be the unit circle bundle, with projection $\pi:X\rightarrow M$. Then $X$ is a contact manifold, 
with contact form given by the connection 1-form  $\alpha$. 
Since $G$ and $T$ preserve the Hermitian metric $h$ on $A^{\vee}$, they act on $X$. Furthermore, 
as both linearized actions preserve the unique compatible connection
$\nabla_{A}$, both actions leave $\alpha$ invariant.

The actions of $G$ and $T$ on $X$ preserve the volume form $dV_{X}=\alpha\wedge \pi^{*}\big(\frac{dV_{M}}{2\pi}\big)$ on $X$, whence they induce commuting unitary representations of $G$ and $T$ on $L^{2}(X)$, which preserve the Hardy space $H(X)=L^2(X)\cap \kernel{(\overline{\partial}_{b})}$.

By virtue of the Peter-Weyl Theorem, we may then unitarily and equivariantly decompose $H(X)$ over the irreducible representations of $G$ and $T$, respectively.  
For every $\nu_{G}\in \widehat{G}$ we define $H(X)_{\nu_{G}}^{G}\subseteq H(X)$ be the maximal subspace equivariantly isomorphic to a direct sum of copies of $V_{\nu_{G}}$. In the same way we define $H(X)_{\nu_{T}}^{T}$. So decomposing the Hardy space of $X$ unitarily and equivariantly over the irreducible representations of $T$ and $G$, we have:

\begin{equation}
\label{pacbell}
H(X)=\bigoplus_{\nu_{T}\in\mathbb{Z}^{d_{T}}}H(X)_{\nu_{T}}^{T}=\bigoplus_{\nu_{G}\in\widehat{G}}H(X)_{\nu_{G}}^{G}.
\end{equation} 

Similarly, under the previous assumptions there is an holomorphic Hamiltonian action of the product $P=G\times T$, and a corresponding unitary representation, so that we also have: 

\begin{equation}
\label{kayak3}
H(X)=\bigoplus_{\nu_{G}\in\widehat{G},\nu_{T}\in\mathbb{Z}^{d_{T}}}H(X)_{\nu_{G},\nu_{T}}^{G\times T},
\end{equation}
\noindent
where $H(X)_{\nu_{G},\nu_{T}}^{G\times T}=H(X)_{\nu_{G}}^{G}\cap H(X)_{\nu_{T}}^{T}$.

Under the assumption that $\mathbf{0}\not\in\Phi_{T}$, we have that $\dimensione(H(X)_{\nu_{T}}^{T})<+\infty$ for each $\nu_{T}\in\mathbb{Z}^{d_{T}}$.

\begin{definition}
Given a pair of irreducible weights $\nu_{G}$ and $\nu_{T}$ for $G$ and $T$, respectively, we shall denote by $\widetilde{\Pi}_{\nu_{G},\nu_{T}}:L^{2}(X)\rightarrow H(X)_{\nu_{G},\nu_{T}}$ the orthogonal projector, and refer to its Schwartz kernel as the level $(\nu_{G},\nu_{T})$-Szeg\"{o} projector of $X$ (with the two actions understood). 
In terms of an orthonormal basis $\left\{s_{j}^{(\nu_{G},\nu_{T})}\right\}_{j=1}^{N_{\nu_{G},\nu_{T}}}$ of $H(X)_{\nu_{G},\nu_{T}}$, it is given by:

\begin{equation}
\label{schwartz}
\widetilde{\Pi}_{\nu_{G},\nu_{T}}(x,y)=\sum_{j}\widehat{s}_{j}^{(\nu_{G},\nu_{T})}(x)\overline{\widehat{s}_{j}^{(\nu_{G},\nu_{T})}(y)}. 
\end{equation}
\end{definition}

In this paper we shall consider the local asymptotics of the equivariant Szeg\"{o} kernels $\widetilde{\Pi}_{\nu_{G},k\nu_{T}}$, where the irreducible representation of $T$ 
tends to infinity along a ray, and the irreducible representation of $G$ is held fixed. To this end, we shall use a combination of the techniques in 
$\cite{primo}$ and $\cite{quarto}$.\\

\begin{observation}
The smooth function $x\mapsto \widetilde{\Pi}_{\nu_{G},k\nu_{T}}(x,x)$ descends to a smooth function on $M$. 
\end{observation}

\begin{remark}
We shall use the following notation for the various equivariant Szeg\"{o} kernels coming into play:

\begin{itemize}
 \item[1] $\widetilde{\Pi}_{\nu_{G},k\nu_{T}}$ in the general case, under the action of $G\times T$;
 \item[2] $\widetilde{\Pi}_{k\nu_{T}}$ in the case of $G$ trivial;
 \item[3] $\Pi_{k}$ in the case of $G$ trivial and $T=S^{1}$ with $\Phi_{T}=1$;
 \item[4] $\Pi_{\nu_{G},k}$ in the case of $T=S^{1}$ and $\Phi_{T}=1$;
 \item[5] $\widetilde{\Pi}_{\nu_{G},k}$ and $\widetilde{\Pi}_{k}$ also for the case $T=S^{1}$ and not necessarily $\Phi_{T}=1$.
\end{itemize}
\end{remark}

\noindent
A key tool used in the proofs are the Heisenberg local coordinates centered at $x\in X$ (see $\cite{terzo}$). 
We denote this system of coordinates by: 

$$ \gamma_{x}:(\theta,\vv)\in (-\pi,+\pi)\times B_{2d_{M}}(\mathbf{0},\delta)\mapsto x+(\theta,\vv)\in X, $$
\noindent
here $B_{2d_{M}}(\mathbf{0},\delta)$ is the open ball of $\mathbb{R}^{2d_{M}}$ centered at the origin with radius $\delta>0$. 
We have that $\theta$ is an angular coordinate along the circle fiber and $\vv\in T_{m}M$ a local coordinate on $M$. We shall also set $x+\vv=x+(0,\vv)$.
Given the choice of HLC centered at $x\in X$, there are induced unitary isomorphisms $T_{x}X\cong \mathbb{R}\oplus \mathbb{R}^{2d_{M}}$ and $T_{m}M \cong \mathbb{R}^{2d_{M}} \cong \mathbb{C}^{d_{M}}$.

Therefore, each equivariant Szeg\"{o} kernel $\widetilde{\Pi}_{\nu_{G},\nu_{T}}$ is a smoothing operator, with $\mathcal{C}^{\infty}$ Schwartz kernel given by $(\ref{schwartz})$.

We shall make the following three transversality assumptions on the moment maps:

\begin{itemize}
\item[1] $\mathbf{0}\not\in \Phi_{T}(M)$ and $\Phi_{T}$ is transversal to the ray $\mathbb{R}_{+}\cdot\nu_{T}$, so that $M_{\nu_{T}}^{T}=\Phi_{T}^{-1}(\mathbb{R}_{+}\cdot\nu_{T})\subseteq M$ is a compact, $G\times T$-invariant and connected submanifold of dimension $2d_{M}+1-d_{T}$. This is equivalent to requiring that the action of $T$ on $X$ be locally free on the inverse of $M^{T}_{\nu_{T}}$ (see $\cite{quarto}$);
\item[2] $\mathbf{0}\in\mathfrak{g}^{\vee}$ is a regular value of $\Phi_{G}$, so that $M_{0}^{G}=\Phi_{G}^{-1}(\mathbf{0})\subseteq M$ is a compact, $G\times T$-invariant and connected submanifold of dimension $2d_{M}-d_{G}$;
\item[3] the two submanifolds $M_{\nu_{T}}^{T}$ and $M_{0}^{G}$ are mutually transversal. 
\end{itemize}

These conditions imply the following (which is what we shall really be using). Since the two actions commute, they give rise to an action of 
the product group $P=G\times T$, which is also holomorphic and Hamiltonian, with moment map:

$$ \Phi_{P}=(\Phi_{G},\Phi_{T}):M\rightarrow \mathfrak{g}^{\vee}\oplus\mathfrak{t}^{\vee}\cong \mathfrak{p}; $$
\noindent
then $\mathbf{0}\not\in \Phi_{P}(M)$, and $\Phi_{P}$ is transversal to the ray $\mathbb{R}_{+}\cdot \big(\mathbf{0},\nu_{T}\big)$. In particular,

$$ M_{0,\nu_{T}}=\Phi^{-1}_{P}\left(\mathbb{R}_{+}\cdot(\mathbf{0},\nu_{T})\right) $$ 

is a smooth connected submanifold of $M$ with codimension $d_{G}+d_{T}-1$.
Leaving aside that here $G$ is not required to be a torus, these hypothesis are similar in nature to the hypothesis in $\cite{quarto}$, applied however to $P$ rather than $T$. Unlike $\cite{quarto}$, where the local scaling asymptotics for representations along a ray $k\nu_{T}$ are considered, we shall study the local scaling asymptotics of doubly equivariant pieces of $\Pi$ associated to pair of representations $(\nu_{G},k\nu_{T})$, where only one of the representations drifts to infinity, while the other is held fixed. Let us also remark that when $T=S^{1}$ and $\Phi_{T}=1$, we are reduced to considering the isotypical components of the spaces of holomorphic global sections $H^{0}(M,A^{\otimes k})$ under the action of $G$, as in $\cite{primo}$ and $\cite{quindicesimo}$.
We will find asymptotic expansions that generalize and combine the previous cases.

Let $N_{m}$ we denote the normal bundle to $\Phi_{P}^{-1}\big(\mathbb{R}_{+}(\mathbf{0},\nu_{T})\big)$ then $N_{m}$ is naturally isomorphic to $J_{m}\big(\kernel{(\Phi_{P}(m))}\big)$ (see $\cite{quarto}$ section 2.2).
The transversality condition is equivalent to require the injectivity of the evaluation map (see as before $\cite{quarto}$). Now for every $m\in M_{0,\nu_{T}}$ we have two Euclidean structures on 

$$\kernel{(\Phi_{P}(m))}=\mathfrak{g}\times\kernel{(\Phi_{T}(m))}\subseteq\mathfrak{g}\oplus\mathfrak{t},$$

\noindent
one induced from $\mathfrak{g}\oplus\mathfrak{t}$ and the second from $T_{m}M$.
Let $D(m)$ the matrix representing the latter Euclidean product on $N_{m}$, with respect to an orthonormal basis. 
Then $D(m)$ is indipendent of the choice of an orthonormal basis for $\mathfrak{g}\times\kernel{(\Phi_{T}(m))}$, 
and it determines a positive smooth function on $M_{0,\nu_{T}}$. As in $\cite{quarto}$ we have the following definition.

\begin{definition}
Define $\mathcal{D}\in\mathcal{C}^{\infty}(M_{0,\nu_{T}})$, with $m\in M_{0,\nu_{T}}$, by setting 

$$ \mathcal{D}(m)=\sqrt{\det{D(m)}}.$$
\end{definition}

We consider $M_{0,\nu_{T}}$ and for the decomposition of the tangent space $T_{m}M$ we have that:

\begin{equation}
\label{labb}
T_{m}M=H_{m}\oplus V_{m}\oplus N_{m}
\end{equation}
\noindent
where, given $J_{m}:T_{m}M\rightarrow T_{m}M$ the complex structure, we have that:  

\begin{equation}
\label{antoniosalierii}
V_{m}=\mathfrak{g}_{M}(m)\oplus\valutation{(\kernel{(\Phi_{T}(m))})}, \ \ N_{m}=J_{m}(V_{m}), \ \  H_{m}=[V_{m}\oplus N_{m}]^{\perp},
\end{equation}
\noindent
are respectively the vertical, the transversal and the horizontal part. Given $m\in M_{0,\nu_{T}}$ and $\vv\in T_{m}M$ we can decompose $\vv$ uniquely as $\vv=\vv_{h}+\vv_{v}+\vv_{t}$ with $\vv_{v}\in V_{m}$, $\vv_{t}\in N_{m}$ and $\vv_{h}\in H_{m}$. 
The scaling asymptotics of the equivariant Szeg\"{o} kernels, that we will see later, are controlled by a quadratic exponent in the components $\vv_h,\vv_v,\vv_t$ of a tangent vector at a given $m\in M_{0,\nu_T}$ (viewed as a small displacement from $m$).

\begin{definition}
Let $x\in X$ and $v_{l}=(\theta_{l},\vv_{l})\in T_{x}X$ with $l=1,2$. We define $H: TX\oplus TX\rightarrow \mathbb{C}$ as

\begin{equation*}
\begin{multlined}[t][12.5cm]
H(v_{1},v_{2})=\\
\lambda_{\nu_{T}}\left(i[\omega_{m}\left(\vv_{1v},\vv_{1t}\right)-\omega_{m}\left(\vv_{2v},\vv_{2t}\right)]+i\omega_{m}\left(\frac{(\theta_{2}-\theta_{1})}{\|\Phi_{T}(m)\|}\eta_{Mh}(m),\vv_{1h}+\vv_{2h}\right)-\right.\\
\left.-i\omega_{m}\left(\vv_{1h},\vv_{2h}\right)-\|\vv_{1t}\|^2-\|\vv_{2t}\|^2-\frac{1}{2}\left\|\vv_{1h}-\frac{(\theta_{2}-\theta_{1})}{\|\Phi_{T}(m)\|}\eta_{Mh}(m)-\vv_{2h}\right\|^2\right)
\end{multlined}
\end{equation*}
\noindent
with $\eta_{Mh}(m)$ the unitary generator of $\kernel{(\Phi_{P}(m))}^{\perp}$ such that $\langle\eta,\Phi_{P}(m)\rangle=\|\Phi_{T}(m)\|$ and $\lambda_{\nu_{T}}=\frac{\|\nu_{T}\|}{\|\Phi_{T}(m)\|}$.
\end{definition}

\begin{theorem}[main Theorem]
\label{teo:secondo}
Under the previous assumptions fix $\nu_{G}\in\widehat{G}$ and consider 
$\nu_{T}\in\mathbb{Z}^{d_{T}}$, assume that $\Phi_{P}$ is transversal to the ray $\mathbb{R}_{+}\cdot(\mathbf{0},\nu_{T})$. We have: 

\begin{itemize}
\item[$1)$] If $C,\delta>0$, and
$$\max{\{\distanza_{M}(\pi(x),M_{0,\nu_{T}}),\distanza_{M}(\pi(y),M_{0,\nu_{T}})\}}\geq Ck^{\delta-\frac{1}{2}},$$ 
\noindent
then $\widetilde{\Pi}_{\nu_{G},k\nu_{T}}=O(k^{-\infty})$.
\item[$2)$] Uniformly in $x\in X_{0,\nu_{T}}$ and $v_{l}\in T_{x}X$ with $\|\vv_{l}\|\leq Ck^{\frac{1}{9}}$, as $k\rightarrow +\infty$ we have: 

\begin{equation*}
\begin{multlined}[t][12.5cm]
\widetilde{\Pi}_{\nu_{G},k\nu_{T}}\left(x+\frac{\vv_{1}}{\sqrt{k}},x+\frac{\vv_{2}}{\sqrt{k}}\right)\sim \frac{1}{(\sqrt{2}\pi)^{d_{T}-1}}d_{\nu_{G}}2^{\frac{d_{G}}{2}}\cdot  \\
\cdot\left(\frac{k}{\pi}\|\nu_{T}\|\right)^{d_{M}-\frac{d_{P}}{2}+\frac{1}{2}}\left(\sum_{j=1}^{N_{x}}\chi_{\nu_{G}}(g_{j}^{-1})e^{-ik\vartheta_{j}\nu_{T}}e^{H(v_{1}^{j},v_{2})}\right)\cdot \frac{e^{-i\sqrt{k}(\theta_{2}-\theta_{1})\lambda_{\nu_{T}}}}{\mathcal{D}(m)}\cdot \\
\cdot\frac{1}{\|\Phi_{T}\|^{d_{M}+1-\frac{d_{P}}{2}+\frac{1}{2}}}\left(1+\sum_{l\geq 1}R_{\nu_{G},l}(m,v_{1}^{j},v_{2})k^{-\frac{l}{2}}\right)
\end{multlined}
\end{equation*}

where $d_{P}=d_{G}+d_{T}$, $v_{1}^{j}$,  $\vv^{j}$ denote the monodromy representation 
$F_{x} \rightarrow GL(T_{m}M)$, such that for every $j=1,\cdots, N_{x}$, $\vv\in T_{m}M$ we have $p_{j} \mapsto d_{m}\mu^{P}_{p_{j}}(\vv)=\vv^{(j)}\in T_{m}M$. Where $F_{x}$ is the stabilizator of $P$ in $x$ and $R_{\nu_{G},l}$ are polynomials in $v_{1}^{j},v_{2}$ with coefficients depending on $x$,$\nu_{G}$ and $\nu_{T}$.
\item[$3)$]More in general, for every $p_{0}\in P$, denoting $P\cdot x$ the orbit of $x\in X_{0,\nu_{T}}$, then the following expansion holds for $k\rightarrow +\infty$: 

\begin{equation}
\label{generalizationcase2}
\begin{multlined}[t][12.5cm]
\widetilde{\Pi}_{\nu_{G},k\nu_{T}}\left(x+\frac{u_{1}}{\sqrt{k}},p_{0}\cdot\left(x+\frac{u_{2}}{\sqrt{k}}\right)\right)\sim \frac{1}{(\sqrt{2}\pi)^{d_{T}-1}}d_{\nu_{G}}2^{\frac{d_{G}}{2}}\cdot\\
\cdot\left(\frac{k}{\pi}\|\nu_{T}\|\right)^{d_{M}-\frac{d_{P}}{2}+\frac{1}{2}}\sum_{j=1}^{N_{x}}\overline{\chi_{\nu_{P}}(p_{j}p^{-1}_{0})}e^{H(\vv^{j}_{1},\vv_{2})}\cdot \frac{e^{-i\sqrt{k}(\theta_{2}-\theta_{1})\lambda_{\nu_{T}}}}{\mathcal{D}(m)}\cdot \\
\cdot\frac{1}{\|\Phi_{T}\|^{d_{M}+1-\frac{d_{P}}{2}+\frac{1}{2}}}\left(1+\sum_{l\geq 1}R_{\nu_{G},l}(m,\vv^{j}_{1},\vv_{2})k^{-\frac{l}{2}}\right),
\end{multlined}
\end{equation}
\noindent
where $p_{j}\in F_{x}$ and $u_{j}=(\theta_{j},\vv_{j})$ for $j=1,2$.
\end{itemize}
\end{theorem}

The previous result describes the asymptotics of $\widetilde{\Pi}_{\nu_{G},k\nu_{T}}$ in a shrinking neighborhood of the orbit $P\cdot x$, where $x\in X_{0,\nu_{T}}$. It is complemented by the following:

\begin{proposition}
\label{prop:complement}
Suppose $x\in X_{0,\nu_{T}}$ and $\varepsilon,D>0$. Then uniformly for $\distanza_{X}(y,P\cdot x)\geq Dk^{\varepsilon-1/2}$ we have 

$$\widetilde{\Pi}_{\nu_{G},k\nu_{T}}(x,y)=O(k^{-\infty}).$$
\end{proposition}

\subsection{Special cases and relation to prior work}

Before continuing our exposition, it is in order to digress on the relation of our results to prior work in this area. Let us focus on the following two special cases:

\begin{itemize}
\item[a)] $T=S^{1}$ acts in the standard manner (with $\Phi_{T}=1$);
\item[b)] $G$ is trivial.
\end{itemize}

Let us first consider the case a), and to fix ideas let us start with the case where $G$ is trivial.
Let $\rho(\cdot,\cdot)$ be a system of Heisenberg local coordinates for $X$ centered at $x$.
We have for $X$ centered at $x$, inducing a unitary isomorphism $(T_{m}M,h_{m})$ with $\mathbb{C}^{d_{M}}$ the standard Hermitian structure.
In Theorem 3.1 of $\cite{terzo}$ and in $\cite{dodicesimo}$, for $\vv_{1},\vv_{2}\in B(\mathbf{0},1)\subseteq \mathbb{C}^{d_{M}}$, $\theta\in (-\pi,\pi)$ and $k\rightarrow +\infty$ the following expansion has been determined for the level $k$ of Szeg\"{o} kernel $\Pi_{k}$ (see also $\cite{tredicesimo}$):

\begin{equation}
\label{standardcaseTYZ55}
\begin{multlined}[t][12.5cm]
\Pi_{k}\left(\rho\left(\theta,\frac{\vv_{1}}{\sqrt{k}}\right),\rho\left(\theta',\frac{\vv_{2}}{\sqrt{k}}\right)\right)\sim\left(\frac{k}{\pi}\right)^{d_{M}}\cdot e^{ik(\theta-\theta')+\psi_{2}(\vv_{1},\vv_{2})}\cdot\\ \cdot\left(1+\sum_{j\geq 1}a_{j}(x,\vv_{1},\vv_{2})k^{-\frac{j}{2}}\right)
\end{multlined}
\end{equation}

\noindent
where 

$$\psi_{2}(\vv_{1},\vv_{2})=\vv_{1}\cdot\overline{\vv_{2}}-\frac{1}{2}(\|\vv_{1}\|^2+\|\vv_{2}\|^2)$$

\noindent
and $a_{j}$ are polynomials in $\vv_{1}$ and $\vv_{2}$.

\begin{observation}
Another way to write $\psi_{2}$ is:

$$ \psi_{2}(\vv_{1},\vv_{2})=-i\omega_{m}\left(\vv_{1},\vv_{2}\right)-\frac{1}{2}\|\vv_{1}-\vv_{2}\|^{2},$$
\noindent
in this form we can see directly the real and the imaginary part of $\psi_{2}$ observing that it is responsible to the exponential decay near the diagonal. 
\end{observation}

Now let us consider the Hamiltonian action of a compact Lie group $G$ on $M$ and suppose that $\mathbf{0}\in\mathfrak{g}^{\vee}$ is a regular value of the moment map. Then $\Pi_{\nu_{G},k}(x,x)$ is rapidly decreasing away from $\Phi_{G}^{-1}(\mathbf{0})$, and assuming $\Phi_{G}(\pi(x))=\mathbf{0}$, under the standard action of $S^{1}$ the following asymptotic expansion holds with $m=\pi(x)$:

\begin{equation}
\label{groooup}
\begin{multlined}[t][12.5cm]
\Pi_{\nu_{G},k}\left(x+\frac{v_{1}}{\sqrt{k}},x+\frac{v_{2}}{\sqrt{k}}\right) \\
\sim \left(\frac{k}{\pi}\right)^{d_{M}-\frac{d_{G}}{2}}e^{\left[Q\left(v_{1v}+v_{1t},v_{2v}+v_{2t}\right)\right]}\sum_{g\in G_{m}}e^{\psi_{2}\left(v_{1gh},v_{2h}\right)}\cdot A_{\nu_{G},k}(g,x)\cdot\\
\cdot\left(1+\sum_{j\geq 1}R_{\nu_{G},j}(m,v_{1g},v_{2})k^{-\frac{j}{2}}\right)
\end{multlined}
\end{equation}
\noindent
where $Q\left(v_{1v}+v_{1t},v_{2v}+v_{2t}\right)=-\|v_{2t}\|^2-\|v_{1t}\|^2+i[\omega_{m}(v_{1v},v_{1t})-\omega_{m}(v_{2v},v_{2t})]$, 
$G_{m}=\{g\in G:\mu_{g}(m)=m\}$, $R_{\nu_{G},j}$ are polynomials in $v_{1},v_{2}$ and 

$$A_{\nu_{G},k}(g,x)=2^{\frac{d_{G}}{2}}\frac{\dimensione(V_{\nu_{G}})}{V_{\eff}{(\pi(x))}}\frac{1}{|G_{\pi(x)}|}\chi_{\nu_{G}}(g)h_{g}^{k},$$
\noindent
where $V_{\eff}{(\pi(x))}$ is the volume of the fiber above $m$ in $\Phi^{-1}_{G}(\mathbf{0})$ (for more details on the effective potentials see $\cite{venticinquesimo}$) and here we have set 

$$ v_{1g}=d_{m}\mu^{G}_{g_{j}^{-1}}(v_{1}) $$
\noindent
with $g_{j}$ in the stabilizator of $G$. Obviously $(\ref{groooup})$ reduces to $(\ref{standardcaseTYZ55})$ for trivial $G$.

Let us consider case b). Thus assume that there is a holomorphic Hamiltonian action of a compact torus $T$, and that the moment map determining the linearization is nowhere zero. To fix ideas, let us first consider the case where $T$ is one-dimensional. If $\xi_{M}$ and $\xi_{X}$ are vector fields on $M$ and on $X$ induced by $\mu^{T}$ and $\widetilde{\mu}^{T}$, we have that in the Heisenberg local coordinates $\xi_{X}(x)=(-\Phi_{T}(m),\xi_{M}(m))$ with $m=\pi(x)$. Let $\xi_{X}(x)^{\perp}\subseteq T_{x}X$ be the orthocomplement of $\xi_{X}(x)$. In view of Theorem 1 of $\cite{quarto}$, again working in a system of HLC centered at $x$ and that $v_{l}=(\theta_{l},\vv_{l})\in T_{x}X\cong \mathbb{R}\times T_{m}M$ satisfying $\vv_{l}\in\xi_{X}(x)^{\perp}$, $\|\vv_{l}\|\leq Ck^{1/9}$, as $k\rightarrow +\infty$ we have:

\begin{equation}
\label{tooooruscase}
\begin{multlined}[t][12.5cm]
\widetilde{\Pi}_{k}\left(x+\frac{\vv_{1}}{\sqrt{k}},x+\frac{\vv_{2}}{\sqrt{k}}\right)\\
\sim \left(\frac{k}{\pi}\right)^{d_{M}}\Phi_{T}(m)^{-(d_{M}+1)}e^{i\sqrt{k}\frac{(\theta_{1}-\theta_{2})}{\Phi_{T}(m)}}\cdot\left(\sum_{t\in T_{m}}t^{k}e^{E(d_{x}\tilde{\mu}^{T}_{t^{-1}}(\vv_{1}),\vv_{2})}\right)\cdot\\
\cdot\left(1+\sum_{j\geq 1}R_{j}(m,v_{1},v_{2})k^{-\frac{j}{2}}\right)
\end{multlined}
\end{equation}
\noindent
for certain smooth functions $R_{j}$, polynomial in the $v_{l}$'s, with 

\begin{equation*}
\begin{multlined}[t][12.5cm]
 E(v_{1},v_{2})=\frac{1}{\Phi_{T}(m)}\left\{i\left[ \frac{(\theta_{2}-\theta_{1})}{\Phi_{T}(m)}\omega_{m}(\xi_{M}(m),\vv_{1}+\vv_{2})-\omega_{m}(\vv_{1},\vv_{2})\right]-\right.\\
 \left.-\frac{1}{2}\left\|\vv_{1}-\vv_{2}-\frac{(\theta_{2}-\theta_{1})}{\Phi_{T}(m)}\xi_{M}(m)\right\|^2\right\}.
\end{multlined}
\end{equation*}

This last result can be generalize to a $d_{T}$-dimensional torus as in Theorem 2 of $\cite{quarto}$. In this Theorem 
we have a result similar to the previous but with the appearance of an additional important invariant, which plays a role
analogous to the effective potential in $(\ref{standardcaseTYZ55})$. 
Suppose that $\Phi_{T}$ is transversal to the ray $\mathbb{R}_{+}\cdot \nu_{T}$. Then the normal space to the inverse image $N_{m}$ at any $m\in M_{\nu_{T}}=\Phi_{T}^{-1}(\mathbb{R}_{+}\cdot \nu_{T})$ is $N_{m}\cong J_{m}\left(\kernel{(\Phi_{T}(m))}\right)$ and the evaluation map $\valutation:\kernel{(\Phi_{T}(m))}\rightarrow T_{m}M$ is injective. Therefore, we have on $\kernel{(\Phi_{T}(m))}$ two Euclidean products, and given two orthonormal basis $\mathcal{B}_{1},\mathcal{B}_{2}$ we can consider the matrix $D(m)$ whose determinant is independent of the choice of the basis. Thus we can let $\mathcal{D}(m)=\sqrt{\det{D(m)}}$.
Considering $\nu_{T}\in \mathbb{Z}^{d_{T}}$, as $k\rightarrow +\infty$ we have:

\begin{equation}
\label{tooooruscase2}
\begin{multlined}[t][12.5cm]
\widetilde{\Pi}_{k\nu_{T}}\left(x+\frac{\vv_{1}}{\sqrt{k}},x+\frac{\vv_{2}}{\sqrt{k}}\right)\\
\sim \left(\frac{1}{(\sqrt{2}\pi)^{d_{T}-1}}\right)\left(\|\nu_{T}\|\frac{k}{\pi}\right)^{d_{M}+\frac{1-d_{T}}{2}}\frac{1}{(\|\Phi_{T}\|)^{d_{M}+1+\frac{1-d_{T}}{2}}\mathcal{D}(m)}\\
e^{i\sqrt{k}\frac{(\theta_{1}-\theta_{2})}{\Phi_{T}(m)}}\cdot\left(\sum_{t\in T_{m}}t^{k}e^{H_{m}(d_{x}\tilde{\mu}^{T}_{t^{-1}}(v_{1}),v_{2})}\right)\cdot\\
\cdot\left(1+\sum_{j\geq 1}R_{j}(m,v_{1},v_{2})k^{-\frac{j}{2}}\right)
\end{multlined}
\end{equation}
\noindent
with 

\begin{equation*}
 H_{m}(v_{1},v_{2})=\frac{\|\nu_{T}\|}{\|\Phi_{T}\|}\left[-i\omega_{m}(\vv_{1},\vv_{2})-\|\vv_{1}\|^2-\|\vv_{2}\|^2\right].
\end{equation*}

In this paper, we shall pair these situations. More precisely, we shall assume given actions of $G$ and $T$ as above, compatible in the sense that they commute, and consider the resulting asymptotics relative to a pair $(\nu_{G}, k\nu_{T})$ of irreducible characters, where $\nu_{G}$ is held fixed, and $k\nu_{T}\rightarrow +\infty$ along an integral ray.

We consider the case of a $d_{T}$-dimensional torus. We have shown at the beginning the fundamental result of this work. Now we present some observations.

\begin{observation}
If $G$ is trivial, $d_{G}=0$ the leading term is:

\begin{equation*}
\begin{split}
\frac{e^{-ik\vartheta_{j}\nu_{T}}}{(\sqrt{2}\pi)^{d_{T}-1}}\left(\frac{\|\nu_{T}\|k}{\pi}\right)^{d_{M}-\frac{d_{T}-1}{2}}\cdot\frac{1}{\|\Phi_{T}(m)\|^{d_{M}+1-\frac{d_{T}}{2}+\frac{1}{2}}\mathcal{D}(m)}e^{\lambda_{\nu_{T}}\left[-i\omega_{m}(\vv_{1}^{j},\vv_{2})-\frac{1}{2}\|\vv_{1}^{j}-\vv_{2}\|^2\right]}
\end{split}
\end{equation*}
\noindent
and we're back to the equation $(\ref{tooooruscase2})$ when $\theta_{1}=\theta_{2}=0$.
\end{observation}

\begin{observation}
In the case $T=S^{1}$ with the standard action we have $\lambda_{\nu_{T}}=1$, 
and when $\theta_{1}=\theta_{2}=0$ the result is the formula $(\ref{groooup})$.
\end{observation}

Now we present a Theorem that is the diagonal version, without scaling of the point $2)$ of the main Theorem.

\begin{theorem}
\label{teo:diagonal}
Under the hypothesis of the main Theorem, for $m=\pi(x)\in M_{0,\nu_{T}}$, as $k\rightarrow +\infty$ we have:

\begin{equation}
\label{principal44}
\begin{multlined}[t][12.5cm]
\widetilde{\Pi}_{\nu_{G},k\nu_{T}}(x,x)\\
\sim\frac{d_{\nu_{G}}2^{d_{G}/2}}{\left(\sqrt{2}\pi\right)^{d_{T}-1}}\cdot\left(\frac{\|\nu_{T}\|k}{\pi}\right)^{d_{M}+\frac{1-d_{P}}{2}}\cdot\sum_{j=1}^{N_{x}}\chi_{\nu_{G}}(g_{j}^{-1})e^{-ik\vartheta_{j}\nu_{T}}\\
\cdot\frac{1}{\mathcal{D}(m)\|\Phi_{T}(m)\|^{d_{M}+1+\frac{1-d_{P}}{2}}}\cdot\left(1+\sum_{l\geq 1}B_{l}(m)k^{-l}\right)
\end{multlined}
\end{equation}
\noindent
with $B_{l}$ that are smooth functions on $M_{0,\nu_{T}}$.
\end{theorem}

\begin{corollary}
\label{cor:quarto}
Under the assumptions of Theorem $~\ref{teo:secondo}$,

\begin{equation*}
\begin{split}
\lim_{k\rightarrow +\infty}\left(\frac{\pi}{\|\nu_{T}\|k}\right)^{d_{M}-d_{P}+1}&\dimensione{\left(H(X)_{\nu_{G},k\nu_{T}}\right)}=\\
=& \frac{d_{\nu_{G}}^{2}}{(2\pi)^{d_{T}-1}}\cdot\int_{M_{0,\nu_{T}}}\frac{\|\Phi_{T}(m)\|^{-(d_{M}+1)+d_{P}-1}}{\mathcal{D}(m)}dV_{M}(m).
\end{split}
\end{equation*}
\end{corollary}

As a very special example, we observe that when $d_{T}=1$ and $T^{1}=S^1$ acts trivially on $M$ with moment map $\Phi_{T}=1$, we have $H(X)_{k}$ the $k$-th isotypical component for the standard $S^1$-action on $X$, which is naturally and unitarily isomorphic to $H^{0}(M,A^{\otimes k})$. 
In this case we have the celebrated Tian-Yau-Zelditch expansion. For this result we refer to the work of Zelditch in $\cite{diciannovesimo}$ (see also $\cite{ventitreesimo}$,  $\cite{ventiquattresimo}$, $\cite{trentunesimo}$, $\cite{trentaseiesimo}$ and, for its near diagonal rescaled generalizations see $\cite{dodicesimo}$, $\cite{
terzo}$, $\cite{tredicesimo}$ and $\cite{quattordicesimo}$). 

\subsection{Applications to Toeplitz operator kernels}

By way of application, motivated by the standard Berezin-Toeplitz quantization of a classical observable (see $\cite{trentesimo}$, $\cite{diciottesimo}$, $\cite{trentaquattresimo}$, $\cite{ventottesimo}$, $\cite{ventinovesimo}$, $\cite{trentatreesimo}$, $\cite{trentottesimo}$, $\cite{trentacinquesimo}$, $\cite{trentaduesimo}$ and $\cite{trentasettesimo}$), let us consider the scaling asymptotics of the equivariant components of certain Toeplitz operators (we will consider Toeplitz operators in the sense of $\cite{undicesimo}$). Given $f\in\mathcal{C}^{\infty}(M)$ and assuming for simplicity that $f$ is invariant under the action of the product group $P=G\times T$, we can consider the Toeplitz operators $T_{\nu_{G},k\nu_{T}}[f]=\widetilde{\Pi}_{\nu_{G},k\nu_{T}}\circ M_{f}\circ \widetilde{\Pi}_{\nu_{G},k\nu_{T}}$, where $M_{f}$ denotes multiplication by $f\circ \pi$. Then $T_{\nu_{G},k\nu_{T}}[f]$ is a self-adjoint endomorphisms of $H(X)_{\nu_{G},k\nu_{T}}$.

Given that $\mathbf{0}\not\in\Phi_{P}$, the equivariant Toeplitz operator $T_{\nu_{G},k\nu_{T}}[f]$ is smoothing, and its distributional kernel is given by the following two alternative expressions:

\begin{equation}
\label{toeplitz}
\begin{multlined}[t][12.5cm]
T_{\nu_{G},k\nu_{T}}[f](x,x')=\int_{X}\widetilde{\Pi}_{\nu_{G},k\nu_{T}}(x,y)f(y)\widetilde{\Pi}_{\nu_{G},k\nu_{T}}(y,x')dV_{X}(y)\\
=\sum_{j}T_{\nu_{G},k\nu_{T}}[f](s_{j}^{k}(x))\overline{(s_{j}^{k}(x'))}
\end{multlined}
\end{equation}
\noindent
with $x,x'\in X_{0,\nu_{T}}$ and $s_{j}^{k}$ an orthonormal basis of $H(X)_{\nu_{G},k\nu_{T}}$.
We will see that $T_{\nu_{G},k\nu_{T}}[f](x,x')$ has asymptotic expansions near the diagonal similar to the one for $\widetilde{\Pi}_{\nu_{G},k\nu_{T}}$.
Note that with $f(y)$ we denote $f(\pi(y))$ and that every $f\in \mathcal{C}^{\infty}(M)$ lifts to an invariant function $f(x)$ on $X$.
For the sake of simplicity, we shall focus on points of the form $(x+n,x+n)$ (with rescaling), as usual, in a system of Heisenberg local coordinates centered at $x$, where $n$ is a normal vector to the $P$-orbit of $x$ and we shall make the extra assumption that the stabilizer of $x$ in $P$ is trivial. Notice that any point sufficiently close to $P\cdot x$ may be written in this manner, possibly replacing $x$ with $p\cdot x$ for some $p\in P$.

\begin{theorem}
\label{teo:app2}
Assume that $\mathbf{0}\not\in \Phi_{P}$, $f\in\mathcal{C}^{\infty}(M_{0,\nu_{T}})$ is $\mu^{P}$-invariant and that the stabilizer of $P$ in $x$ is trivial. Suppose $x\in X_{0,\nu_{T}}$ and fix a system of HLC centered at $x$. Let $m=\pi(x)$. Then we have:

\begin{itemize}
\item[$1)$]If $C,\delta>0$ and 
$$\max{\{\distanza_{M}(\pi(x),M_{0,\nu_{T}}),\distanza_{M}(\pi(y),M_{0,\nu_{T}})\}}\geq Ck^{\delta-\frac{1}{2}},$$ 
\noindent
then $T_{\nu_{G},k\nu_{T}}[f](x,x')=O(k^{-\infty})$.
\item[$2)$]Uniformly in $n_{1}\in N_{x}^{P}=T_{x}(P\cdot x)^{\perp}$ as $k\rightarrow +\infty$:

\begin{equation}
\label{Toeplitz7}
\begin{multlined}[t][12.5cm]
T_{\nu_{G},k\nu_{T}}[f]\left(x+\frac{n_{1}}{\sqrt{k}},x+\frac{n_{1}}{\sqrt{k}}\right)\\
\sim \frac{1}{(\sqrt{2}\pi)^{d_{T}-1}}d_{\nu_{G}}2^{\frac{d_{G}}{2}}\left(\frac{k}{\pi}\|\nu_{T}\|\right)^{d_{M}-\frac{d_{P}}{2}+\frac{1}{2}}f(m)e^{-2\lambda_{\nu_{T}}\|\ttt_{1}\|^2}\cdot\\
\cdot \frac{1}{\mathcal{D}(m)}\cdot\frac{1}{\|\Phi_{T}\|^{d_{M}+1-\frac{d_{P}}{2}+\frac{1}{2}}}\left(\sum_{l\geq 0}k^{-\frac{l}{2}}R_{l}(n_{1},m)\right)
\end{multlined}
\end{equation}
\noindent
with $R_{l}(n_{1},m)$ a polynomial in $n_{1}$ and $\ttt_{1}\in N_{m}=J_{m}(\valutation_{m}(\kernel{(\Phi_{P}(m))}))$.
\end{itemize}
\end{theorem}

\begin{corollary}
\label{cor:quinto}
Under the assumptions of Theorem $~\ref{teo:app2}$,

\begin{equation*}
\begin{split}
\lim_{k\rightarrow +\infty}\left(\frac{\pi}{\|\nu_{T}\|k}\right)^{d_{M}-d_{P}+1}&\mathfrak{T}\left(T_{\nu_{G},k\nu_{T}}[f]\right)=\\
=& \frac{d_{\nu_{G}}^{2}}{(2\pi)^{d_{T}-1}}\cdot\int_{X_{0,\nu_{T}}}\frac{f(\pi(x))\|\Phi_{T}(\pi(x))\|^{-(d_{M}+2-d_{P})}}{\mathcal{D}(\pi(x))}dV_{X}(x),
\end{split}
\end{equation*}
\noindent
where $\mathfrak{T}\left(T_{\nu_{G},k\nu_{T}}[f]\right)$ is the trace of the Toeplitz operator.
\end{corollary}

\section{Examples}

The main Theorem predicts that the diagonal restriction $\widetilde{\Pi}_{\nu_{G},k\nu_{T}}(x,x)$ of the equivariant Szeg\"{o} kernel (which descends to a function on $M$) is rapidly decreasing away from the locus $M_{0,\nu_{T}}$, and grows like $k^{d_{M}+\frac{1-d_{P}}{2}}$ there. Let us illustrate this explicitly by two examples (cfr the computations in $\cite{primo}$). Recall from $\cite{dodicesimo}$ that for $k=1,2,\cdots$ an orthonormal basis of 
$H^{0}(\mathbb{P}^{n},\mathcal{O}_{\mathbb{P}^{n}}(k))$ is $\{s_{J}^{k}\}_{|J|=k}$, where:

\begin{equation}
s_{J}^{k}=\sqrt{\frac{(k+n)!}{\pi^{n}J!}}z^{J}
\end{equation}
\noindent
and where $J!=\prod_{l=0}^{n}j_{l}!$, $z^{J}=\prod_{l=0}^{n}z_{l}^{j_{l}}$.
In the next example we consider a particular product action and we show that outside of $M_{0,\nu_{T}}$ we have the exponential decay of the Szeg\"{o} kernel.

\begin{example}
Let us make $M=\mathbb{P}^{1}$. Let us consider the action of $G=T^{1}$ on $M$ induced by the representation on $\mathbb{C}^{2}$ given by $\mu^{G}(z_{0},z_{1})=w\cdot(z_{0},z_{1})=(w^{-1}z_{0},wz_{1})$, and the action of $T=T^{1}$ induced by the representation given by $\mu^{T}(z_{0},z_{1})=(s^{-1}z_{0},s^{-2}z_{1})$. These actions are holomorphic and Hamiltonian, with moment maps:
 
$$\Phi_{G}(z_{0},z_{1})=\frac{|z_{0}|^2-|z_{1}|^2}{|z_{0}|^2+|z_{1}|^2}$$ 
\noindent
and 

$$\Phi_{T}(z_{0},z_{1})=\frac{|z_{0}|^2+2|z_{1}|^2}{|z_{0}|^2+|z_{1}|^2}.$$

Then we have:

$$\Phi_{G}^{-1}(0)=\{[z_{0}:z_{1}]:|z_{0}|=|z_{1}|\}$$
\noindent
and placing $X=S^{3}\subseteq \mathbb{C}^{2}$ we have 

$$X_{0}=\pi^{-1}\left(\Phi_{G}^{-1}(0)\right)=\left\{(z_{0},z_{1}):|z_{0}|=|z_{1}|=\frac{1}{\sqrt{2}}\right\}\cong S^{1}\times S^{1}$$
\noindent
with a free action of $S^{1}$ on $X_{0}$.  
We have $\nu_{T}=1\in\mathbb{Z}$, $\Phi_{P}^{-1}(\mathbb{R}_{+}\cdot(0,1))=\Phi_{G}^{-1}(0)=\{(z_{0},z_{1}):|z_{0}|=|z_{1}|\}$ and the action of $P$ is given by:

$$\mu^{P}([z_{0}:z_{1}])=(w,s)\cdot (z_{0},z_{1})=\left((ws)^{-1}z_{0},ws^{-2}z_{1}\right).$$

If $|z_{0}|=|z_{1}|$ ($\not=0$) and $(w,s)\cdot (z_{0},z_{1})=(z_{0},z_{1})\Rightarrow ws=1, ws^{-2}=1 \Rightarrow s=s^{-2}$ so $s=e^{\frac{2}{3}\pi ji}$ with $j=0,1,2$ and $w=\frac{1}{s}$ then the action is locally free. We are in the hypothesis of the main Theorem.
We have $s\cdot(z_{0}^{a}z_{1}^{b})=(sz_{0})^{a}(s^{2}z_{1})^{b}=s^{a+2b}z_{0}^{a}z_{1}^{b}$ and then 

$$\widetilde{H}^{T}(X)_{k}=\spaned{\{z_{0}^{a}z_{1}^{b}: a+2b=k\}}.$$

In the other side we have $w\cdot(z_{0}^{a}z_{1}^{b})=(wz_{0})^{a}(w^{-1}z_{1})^{b}=w^{a-b}z_{0}^{a}z_{1}^{b}$ and then 

$$\widetilde{H}^{G}(X)_{\nu_{G}}=\spaned{\{z_{0}^{a}z_{1}^{b}: a=b+\nu_{G}\}}.$$

Thus

$$ \widetilde{H}^{P}(X)_{\nu_{G},k}=\spaned{\left\{z_{0}^{a}z_{1}^{b}: a=b+\nu_{G},a+2b=k\right\}}$$
\noindent
then $a+2b=k \Rightarrow b+\nu_{G}+2b=k \Rightarrow 3b=k-\nu_{G}$ and

$$ \dimensione{\left(\widetilde{H}^{P}(X)_{\nu_{G},k}\right)}=\begin{cases} 0  &\mbox{if} \ \  k\equiv\nu_{G} \mod{3}\\ 
1 & \mbox{if} \ \ k\not\equiv \nu_{G}\mod{3}  \end{cases}. $$

If $k=\nu_{G}+3b$ we have:

$$ \widetilde{H}^{P}_{\nu_{G},\nu_{G}+3b}(X)=\spaned{\left\{z_{0}^{b+\nu_{G}}z_{1}^{b}\right\}}, $$
\noindent
the corresponding Szeg\"{o} projector is:

\begin{equation}
\label{loopspaces}
\begin{multlined}[t][12.5cm]
\widetilde{\Pi}^{P}_{\nu_{G},\nu_{G}+3b}\left((z_{0},z_{1}),(u_{0},u_{1})\right)=\frac{(2b+\nu_{G}+1)!}{\pi(b+\nu_{G})!b!}(z_{0}\overline{u}_{0})^{b+\nu_{G}}(z_{1}\overline{u}_{1})^{b}.
\end{multlined}
\end{equation}

Now consider $z_{j}=u_{j}$ with $|z_{0}|^2+|z_{1}|^2=1$ and we set $x=|z_{0}|^{2}$, $y=y_{b}=\frac{b+\nu_{G}}{2b+\nu_{G}}\rightarrow \frac{1}{2}$ as $b\rightarrow +\infty$.
Using Stirling approximation:

$$ n!\sim \sqrt{2\pi n}\frac{n^{n}}{e^{n}}$$
\noindent
and the projector:

$$
\widetilde{\Pi}^{P}_{\nu_{G},\nu_{G}+3b}\left((z_{0},z_{1}),(z_{0},z_{1})\right)=\frac{(2b+\nu_{G}+1)!}{\pi(b+\nu_{G})!b!}|z_{0}|^{2(b+\nu_{G})}|z_{1}|^{2b}
$$
\noindent
we can find the following asymptotic for the coefficient:

\begin{equation}
\label{loopspaces2}
\begin{multlined}[t][12.5cm]
\frac{(2b+\nu_{G}+1)!}{\pi(b+\nu_{G})!b!}\sim \frac{2}{\sqrt{\pi}}\sqrt{b}\left(\frac{1}{y_{b}}\right)^{b+\nu_{G}}\left(\frac{1}{1-y_{b}}\right)^{b}
\end{multlined}
\end{equation}
\noindent
and for the projector:

\begin{equation}
\label{loopspaces3}
\begin{multlined}[t][12.5cm]
\widetilde{\Pi}^{P}_{\nu_{G},\nu_{G}+3b}\left((z_{0},z_{1}),(z_{0},z_{1})\right)\sim 2\sqrt{\frac{b}{\pi}}\left(\frac{x}{y}\right)^{\nu_{G}}e^{bF(x,y)},
\end{multlined}
\end{equation}
\noindent
where we set $F(x,y)=\log{x}+\log{(1-x)}-\log{y}-\log{(1-y)}=f(x)-f(y)$ with $f(t)=\log{t}+\log{(1-t)}$ and $0<t<1$.
We observe that for $t\rightarrow 0^{+},1^{-}$ we obtain $f(t)\rightarrow -\infty$ and that the derivative:

$$ f'(t)=\frac{1}{t}-\frac{1}{1-t}=0 \Leftrightarrow 1-t= t \Leftrightarrow t=\frac{1}{2} $$
\noindent
with $f(1/2)=-\log{4}$. Thus for $b\gg 0$ and $y=y_{b}$ we have $f(y)=-\log{4}-\delta(b)$ with $\delta(b)>0$ and $\delta(b)\rightarrow 0$ as $b\rightarrow +\infty$. If $x \not=\frac{1}{2}$ we have $f(x)=-\log{4}-\delta(x)$ ( with $\delta(x)>0$ fixed).
Then $F(x,y_{b})=-\delta(x)+\delta(b)\leq -\frac{\delta(x)}{2}$.

Now 

\begin{equation}
\label{loopspaces44}
\begin{multlined}[t][12.5cm]
\left|\widetilde{\Pi}^{P}_{\nu_{G},\nu_{G}+3b}\left((z_{0},z_{1}),(z_{0},z_{1})\right)\right|\sim 2\sqrt{\frac{b}{\pi}}\left(\frac{x}{y}\right)^{\nu_{G}}e^{bF(x,y)}\\
\leq 2\sqrt{\frac{b}{\pi}}\left(\frac{x}{y}\right)^{\nu_{G}}e^{-b\frac{\delta(x)}{2}}=O(b^{-\infty})
\end{multlined}
\end{equation}
\noindent
but $x=\frac{1}{2}$ that is $|z_{0}|=|z_{1}|=\frac{1}{\sqrt{2}}$ and we have:

\begin{equation}
\label{loopspaces55}
\begin{multlined}[t][12.5cm]
\left|\widetilde{\Pi}^{P}_{\nu_{G},\nu_{G}+3b}\left((z_{0},z_{1}),(z_{0},z_{1})\right)\right|\sim 2\sqrt{\frac{b}{\pi}}\left(\frac{x}{y}\right)^{\nu_{G}}e^{bF(1/2,y_{b})}\\
\leq 2\sqrt{\frac{b}{\pi}}\left(\frac{1}{2y_{b}}\right)^{\nu_{G}}e^{b\delta(b)}
\end{multlined}
\end{equation}
\noindent
and considering that we have for $y_{b}$:

\begin{equation}
\label{loopspaces56}
\begin{multlined}[t][12.5cm]
y_{b}=\frac{b+\nu_{G}}{2b+\nu_{G}}=\frac{1+\frac{\nu_{G}}{b}}{2\left(1+\frac{\nu_{G}}{2b}\right)}=\frac{1}{2}+\frac{\nu_{G}}{4b}+O\left(\frac{1}{b^{2}}\right),
\end{multlined}
\end{equation}
\noindent
then $f(y_{b})=-\log{4}+O\left(\frac{1}{b^{2}}\right)$ (because $f'(1/2)=0$) and so follows that $b\delta(b)=O\left(\frac{1}{b}\right)\rightarrow 0$ as $b\rightarrow +\infty$. Thus

\begin{equation}
\label{loopspaces5566}
\begin{multlined}[t][12.5cm]
\left|\widetilde{\Pi}^{P}_{\nu_{G},\nu_{G}+3b}\left((z_{0},z_{1}),(z_{0},z_{1})\right)\right|\sim  2\sqrt{\frac{b}{\pi}}.
\end{multlined}
\end{equation}
\end{example}

Another possible variation similar to the previous is the following. 

\begin{example}
Let us make $M=\mathbb{P}^{2}$. Let us consider the actions of $G=T^{2}$ on $M$ induced by the representation on $\mathbb{C}^{3}$ given by $\mu^{G}(z_{0},z_{1},z_{2})=(w_{1}^{-1}z_{0},w_{1}w_{2}^{-1}z_{1},w_{2}z_{2})$, and the action of $T=T^{1}$ induced by the representation given by $\mu^{T}(z_{0},z_{1},z_{2})=(s^{-1}z_{0},s^{-2}z_{1},s^{-3}z_{2})$. These actions are holomorphic and Hamiltonian, with moment maps:
 
$$\Phi_{G}(z_{0},z_{1},z_{2})=\left(\frac{|z_{0}|^2-|z_{1}|^2}{|z_{0}|^2+|z_{1}|^2},\frac{|z_{1}|^2-|z_{2}|^2}{|z_{1}|^2+|z_{2}|^2}\right)$$ 
\noindent
and 

$$\Phi_{T}(z_{0},z_{1},z_{2})=\frac{|z_{0}|^2+2|z_{1}|^2+3|z_{2}|^{2}}{|z_{0}|^2+|z_{1}|^2+|z_{2}|^{2}}.$$

Then

$$\Phi_{G}^{-1}(\mathbf{0})=\{[z_{0}:z_{1}:z_{2}]:|z_{0}|=|z_{1}|=|z_{2}|\}$$
\noindent
and placing $X=S^{5}\subseteq \mathbb{C}^{3}$ we have:

$$X_{0}=\pi^{-1}\left(\Phi_{G}^{-1}(\mathbf{0})\right)=\left\{(z_{0},z_{1},z_{2}):|z_{0}|=|z_{1}|=|z_{2}|=\frac{1}{\sqrt{3}}\right\}\cong S^{1}\times S^{1}\times S^{1}$$
\noindent
with a free action of $G$ on $X_{0}$.
We have $\nu_{T}=1\in\mathbb{Z}$ and 
$\Phi_{P}^{-1}(\mathbb{R}_{+}\cdot(\mathbf{0},1))=\Phi_{G}^{-1}(\mathbf{0})=\{(z_{0},z_{1},z_{2}):|z_{0}|=|z_{1}|=|z_{2}|\}$ and the action of $P$ is given by:

$$\mu^{P}([z_{0}:z_{1}:z_{2}])=(w,s)\cdot (z_{0},z_{1},z_{2})=\left((w_{1}s)^{-1}z_{0},w_{1}w_{2}^{-1}s^{-2}z_{1},w_{2}s^{-3}z_{2}\right).$$

If $|z_{0}|=|z_{1}|=|z_{2}|$ ($\not=0$) and $(w,s)\cdot (z_{0},z_{1},z_{2})=(z_{0},z_{1},z_{2})\Rightarrow w_{1}s=1, w_{1}w_{2}^{-1}s^{-2}=1, w_{2}s^{-3}=1 \Rightarrow s^{6}=1$ so $s=e^{\frac{2}{6}\pi ji}$ with $j=0,1,2,4,5$ and $w_{1}=\frac{1}{s},w_{2}=\frac{1}{s^{3}}$ then the action is locally free. The hypothesis of the main Theorem are satisfied.
We have $s\cdot(z_{0}^{a}z_{1}^{b}z_{2}^{c})=(sz_{0})^{a}(s^{2}z_{1})^{b}(s^{3}z_{2})^{c}=s^{a+2b+3c}z_{0}^{a}z_{1}^{b}z_{2}^{c}$ and then 

$$\widetilde{H}^{T}(X)_{k}=\spaned{\{z_{0}^{a}z_{1}^{b}z_{2}^{c}: a+2b+3c=k\}}.$$

In analogue way we obtain that 

$$\widetilde{H}^{G}(X)_{\nu_{G}}=\spaned{\{z_{0}^{a}z_{1}^{b}z_{2}^{c}: (a-b,b-c)=(\nu_{1},\nu_{2})=\nu_{G}\}}.$$

Thus

$$\widetilde{H}^{P}(X)_{\nu_{G},k}=\spaned{\left\{z_{0}^{a}z_{1}^{b}z_{2}^{c}: (a-b,b-c)=(\nu_{1},\nu_{2}),a+2b+3c=k\right\}}$$
\noindent
and 

$$\dimensione{\left(\widetilde{H}^{P}(X)_{\nu_{G},k}\right)}=\begin{cases} 0  &\mbox{if } k\equiv\nu_{1} \mod{6}\\ 
1 & \mbox{if } k\not\equiv \nu_{1}\mod{6}  \end{cases}.$$

If $k=6c+\nu_{1}+3\nu_{2}$ we have:

$$\widetilde{H}^{P}(X)_{(\nu_{1},\nu_{2}),\nu_{1}+3\nu_{2}+6c}=\spaned{\left\{z_{0}^{c+\nu_{1}+\nu_{2}}z_{1}^{c+\nu_{2}}z_{2}^{c}\right\}}.$$

Now consider $z_{j}=u_{j}$ with $|z_{0}|^2+|z_{1}|^2+|z_{2}|^2=1$ and we set $x=|z_{0}|^{2}$, $y=|z_{1}|^2, \ \ z=|z_{2}|^2=1-x-y$.
As before, using Stirling approximation for the projector:

\begin{equation}
\label{loopspaces88}
\begin{multlined}[t][12.5cm]
\widetilde{\Pi}^{P}_{\nu_{G},6c+\nu_{1}+3\nu_{2}}\left((z_{0},z_{1},z_{2}),(z_{0},z_{1},z_{2})\right)\\
=\frac{(3c+\nu_{1}+2\nu_{2}+2)!}{\pi^2(c+\nu_{2}+\nu_{1})!(c+\nu_{2})!c!}(|z_{0}|)^{2(c+\nu_{2}+\nu_{1})}(|z_{1}|)^{2(c+\nu_{2})}(|z_{2}|)^{2c},
\end{multlined}
\end{equation}
\noindent
the following asymptotic holds:

\begin{equation}
\label{loopspaces2acca}
\begin{multlined}[t][12.5cm]
\frac{(3c+\nu_{1}+2\nu_{2}+2)!}{\pi^2(c+\nu_{2}+\nu_{1})!(c+\nu_{2})!c!}\\
\sim \frac{9\sqrt{3}c}{2\pi^3}\left[\frac{3c+\nu_{1}+2\nu_{2}}{c+\nu_{2}+\nu_{1}}\right]^{c+\nu_{2}+\nu_{1}}\left[\frac{3c+\nu_{1}+2\nu_{2}}{c+\nu_{2}}\right]^{c+\nu_{2}}\left[\frac{3c+\nu_{1}+2\nu_{2}}{c}\right]^c
\end{multlined}
\end{equation}
\noindent
and, for the projector:

\begin{equation}
\label{loopspaces3acca}
\begin{multlined}[t][12.5cm]
\widetilde{\Pi}^{P}_{\nu_{G},6c+\nu_{1}+3\nu_{2}}\left((z_{0},z_{1},z_{2}),(z_{0},z_{1},z_{2})\right)\sim \frac{9\sqrt{3}c}{2\pi^3}x^{\nu_{1}+\nu_{2}}y^{\nu_{2}}e^{cF(x,y)},
\end{multlined}
\end{equation}
\noindent
with $F(x,y)= \log{x}+\log{y}+\log{(1-x-y)}$, $0<x,y$ and $x+y<1$.
Studing the partial derivatives we find that we have a critical point $x=y=\frac{1}{3}$ with $F\left(\frac{1}{3},\frac{1}{3}\right)=-\log{27}$.
So we have $F\left(x,y\right)=-\log{27}-\delta(z)$ with $\delta(z)>0$. If $x=y\not=\frac{1}{3}$ we have:

\begin{equation}
\label{loopspaces3141}
\begin{multlined}[t][12.5cm]
\left|\widetilde{\Pi}^{P}_{\nu_{G},6c+\nu_{1}+3\nu_{2}}\left((z_{0},z_{1},z_{2}),(z_{0},z_{1},z_{2})\right)\right|\\
\leq \frac{9\sqrt{3}c}{2\pi^3}x^{\nu_{1}+\nu_{2}}y^{\nu_{2}}e^{-c\delta(z)}=O(c^{-\infty})
\end{multlined}
\end{equation}
\noindent
but $x=y=\frac{1}{3}$, $\delta(z)=0$ and

\begin{equation}
\label{loopspaces3141}
\begin{multlined}[t][12.5cm]
\left|\widetilde{\Pi}^{P}_{\nu_{G},6c+\nu_{1}+3\nu_{2}}\left((z_{0},z_{1},z_{2}),(z_{0},z_{1},z_{2})\right)\right|
\sim\frac{9\sqrt{3}c}{2\pi^3}\left(\frac{1}{3}\right)^{\nu_{1}+\nu_{2}}\left(\frac{1}{3}\right)^{\nu_{2}}.
\end{multlined}
\end{equation}
\end{example}

\section{Preliminaries}

It is know that if $G$ and $T$ both act on a symplectic manifold $M$ in an Hamiltonian fashion with moment maps $\Phi_{G}$ and $\Phi_{T}$ and these actions commute, 
then $P=G\times T$ act on $M$ and the moment map is $\Phi_{P}=\Phi_{G}\oplus\Phi_{T}:M\rightarrow \mathfrak{g}^{\vee}\oplus\mathfrak{t}^{\vee}$.
We give an explicit expression for $H(X)_{\nu_{G},\nu_{T}}$ as

$$H(X)_{\nu_{G},\nu_{T}}=\left\{s\in H(X)_{\nu_{G}}: s(\tilde{\mu}^{T}_{t^{-1}}(x))=t^{\nu_{T}}s(x), \forall x\in X,\forall t\in T^{d_{T}},\forall \nu_{G}\in\widehat{G}\right\}.$$

\subsection{The geometric setting}

We remember that the matrix $D(m)$ represents the Euclidean product on $N_{m}$ with an orthonormal basis. 
It determines a positive smooth function $\mathcal{D}$ on $M_{0,\nu_{T}}$ defined above.

\begin{remark}
Note that there is a relation between the $D(m)$ matrix and the $C(m)$ matrix used in Theorem $\ref{teo:secondo}$. We have in fact that $D(m)=C(m)^{t}\cdot C(m)$ and so 
$\mathcal{D}(m)=\sqrt{\det{D(m)}}=|\det{C(m)}|$.
\end{remark}

Let us consider the symplectic cone $\Sigma \subseteq TX^{\vee}\setminus\{0\}$ sprayed by the connection form $\alpha$:

$$ \Sigma=\{ (x,r\alpha_{x}): x\in X, r>0\}.$$

This cone is important for the microlocal description of Szeg\"{o} kernel (as in $\cite{secondo}$) and in the theory of Toeplitz operators (see $\cite{undicesimo}$). 
We have that the wave front set of 
$\Pi$ is the anti-diagonal:

$$ \Sigma^{\#}=\{(x,r\alpha_{x},x,-r\alpha_{x}): x\in X , r>0\}. $$

Notice that $\Sigma \cong X\times \mathbb{R}_{+}$ in a natural manner.
Let $\omega_{\Sigma}$ the restriction to $\Sigma$ of the symplectic structure on  $TX^{\vee}$. Let $r$ be the cone coordinate on $\Sigma$ and 
$\theta$ be the circle coordinate on $X$, locally defined, and pulled-back to $\Sigma$. Then $\omega_{\Sigma}=d\lambda=dr\wedge \alpha +2r\omega$, 
with $\lambda=r\alpha$. Let $\widetilde{\xi}_{f}$ be the contact lift to $X$ of the Hamiltonian vector field $\xi_{f}$ on $(M,2\omega)$. Then the cotangent flow restricted to $\Sigma$ is generated by $(\widetilde{\xi}_{f},0)$. Thus the cotangent flow on $\Sigma\cong X\times \mathbb{R}$ is $\phi^{\Sigma}_\tau=\phi^{X}_\tau\times id_\mathbb{R}$. It is follows that if $f$ and $g$ Poisson commute on $M$, then their flows on $M$, $X$ and $\Sigma$ also commute, and conversely.
About the product action (referring to $\cite{quarto}$), since by assumption $\Phi_{P}$ is transverse to $\mathbb{R}_{+}\cdot (\mathbf{0},\nu_{T})$, the action of $P$ on $X_{0,\nu_{T}}$ is locally free. We have also that $X_{0,\nu_{T}}$ is invariant for $G\times T$. In fact we have that $\mathbb{R}_{+}\cdot(\mathbf{0},\nu_{T})$ is invariant for the coadjoint action and, with the fact that $\Phi_{P}$ is equivariant, we conclude.

\subsection{Adapted Heisenberg local coordinates}

A key tool used in the proofs are the Heisenberg local coordinates centered at $x\in X$ defined in $\cite{terzo}$.
We choose an adapted holomorphic coordinate system $(z_{1},\cdots , z_{d_{M}})$ for $M$ centered at $\pi(x)$ so $\omega$ 
espressed in $z_{i}$'s at $\pi(x)$ is the standard symplectic structure on $\mathbb{C}^{d_{M}}$ that is 
$\omega\left(\pi(x)\right)=\frac{i}{2}\sum_{j=1}^{d_{M}}dz_{j}\wedge d\overline{z}_{j}$. By the choice of the adapted holomorphic coordinate
system we have the unitary isomorphism $T_{\pi(x)}M\cong\mathbb{C}^{d_{M}}$. Now we choose a preferred local frame $\sigma_{A}$ for $A$ at $\pi(x)$, in the sense of 
$\cite{terzo}$.

The problem now is to find the espression of Heisenberg coordinates under the action of $P=G\times T$. 
The expression of $\tilde{\mu}_{g^{-1}}^{G}\circ\tilde{\mu}_{-\frac{\vartheta}{\sqrt{k}}}^{T}\left(x_{1k}\right)$ in Heisenberg coordinates, where
$x_{1k}=x+\frac{\vv_{1}}{\sqrt{k}}$, $g=\exp_{G}\frac{\varsigma}{\sqrt{k}}$ 
is an element of the group $G$, $e^{i\vartheta}$ an element of the torus $T$ and $\exp_{G}$ is the exponential map, are  given by the following lemma:

\begin{lemma}
\label{lemma:hei}
Suppose $x\in X$, $\Phi_{G}\circ\pi(x)=\mathbf{0}$, and fix a system of Heisenberg local coordinates centered at $x$. Suppose in addition that
$F_{x}=\{(g_{j},t_{j}):j=1,\cdots, N_{x}\}$ is not trivial. 
Then there exist $\mathcal{C}^{\infty}$ functions $\widetilde{B}_{3},\widetilde{B}_{2}:\mathbb{R}^{d_{T}}\times\mathbb{C}^{d_{M}}\times \mathbb{R}^{d_{G}} \rightarrow \mathbb{C}^{d_{M}}$,
vanishing at the origin to third and second order, respectively, such that the following holds.
For $\vv_{1}-\vartheta\xi_{M}(m)\in T_{\pi(x)}M$, where $\vartheta=(\vartheta_{1},\cdots,\vartheta_{d_{T}})$ and $\xi_{M}(m)=\valutation_{m}{\xi}$ 
with $\xi\in \mathfrak{t}$ and $\varsigma_{M}(m)$ the valutation of $\varsigma\in\mathfrak{g}$. As $k \rightarrow +\infty$ the Heisenberg local coordinates of:

$$ \tilde{\mu}^{G}_{e^{-\frac{\varsigma}{\sqrt{k}}}}\circ\tilde{\mu}^{T}_{-\frac{\vartheta}{\sqrt{k}}}\left(\tilde{\mu}^{G}_{g_{j}^{-1}}\circ \tilde{\mu}^{T}_{t_{j}^{-1}}(x_{1k})\right) $$
\noindent
are given by
{ \footnotesize
\begin{equation}
\label{autoscuolavisconti}
\begin{multlined}[t][12.5cm]
\left(\frac{1}{\sqrt{k}}\left(\vartheta\Phi_{T}(m)+\theta_{1}\right)+\frac{1}{k}\omega_{m}\left(\vartheta\xi_{M}(m),\vv_{1}^{j}\right)+\frac{1}{k}\omega_{m}\left(\varsigma_{M}(m),\vv_{1}^{j}\right)+\widetilde{B}_{3}\left(\frac{\vartheta}{\sqrt{k}},\frac{\vv}{\sqrt{k}},\frac{\varsigma}{\sqrt{k}}\right),\right.\\
,\left.\frac{1}{\sqrt{k}}\left(\vv_{1}^{j}-\vartheta\xi_{M}(m)-\varsigma_{M}(m)\right)+\widetilde{B}_{2}\left(\frac{\vartheta}{\sqrt{k}},\frac{\vv}{\sqrt{k}},\frac{\varsigma}{\sqrt{k}}\right)\right).
\end{multlined}
\end{equation}

}\end{lemma}

\emph{Proof.} We apply corollary 2.2 of $\cite{quarto}$ and we have $(\ref{autoscuolavisconti})$.
\hfill $\Box$\\

\begin{remark}
Note that here $\vartheta=(\vartheta_{1},\cdots,\vartheta_{d_{T}})$ with $-\pi<\vartheta_{i}<\pi$ and $\xi=\left.\frac{\partial}{\partial \vartheta}\right|_{0}$. We have: 

$$ \Phi_{l}=\langle\Phi,\xi_{l}\rangle, \ \ \vartheta\cdot\Phi=\sum_{l=1}^{d_{T}}\vartheta_{l}\Phi_{l},  
\ \ \vartheta\cdot\xi=\sum_{l=1}^{d_{T}}\vartheta_{l}\xi_{l}.$$
\end{remark}

\section{Proof of the main Theorem $\bf{\ref{teo:secondo}}$}

$\Proof.$

Proof of $1)$. We consider $(\rho_{\nu_{G}},V_{\nu_{G}})$ an unitary irreducible representation of $G$ and 
we define $\rho_{\nu_{G},k\nu_{T}}:G\times T\rightarrow GL(V_{\nu_{G}})$ as $\rho_{\nu_{G},k\nu_{T}}(g,t)=t^{k\nu_{T}}\rho_{\nu_{G}}(g)$.
We have that $(\rho_{\nu_{G},k\nu_{T}},V_{\nu_{G}})$ is an unitary irreducible representation of $G\times T$ with character 
$\chi_{\nu_{G},k\nu_{T}}(g,t)=\tr{(\rho_{\nu_{G},k\nu_{T}}(g,t))}=t^{k\nu_{T}}\tr{(\rho_{\nu_{G}}(g))}=t^{k\nu_{T}}\chi_{\nu_{G}}(g)$.

Assuming that $\mathbf{0}\not\in \Phi_{T}(M)$ we have that $H(X)_{k\nu_{T}}$ is finite dimensional, then $H(X)_{\nu_{G},k\nu_{T}}\subseteq H(X)_{k\nu_{T}}$
and $\widetilde{\Pi}_{\nu_{G},k\nu_{T}}\in\mathcal{C}^{\infty}(X\times X)$. We want study the asymptotic behavior of $\widetilde{\Pi}_{\nu_{G},k\nu_{T}}$ with $k\rightarrow +\infty$. 
Since $\widetilde{\Pi}_{\nu_{G},k\nu_{T}}$ is the the composition of $\Pi:L^2(X)\rightarrow H(X)$ and the orthogonal projector of $H(X)$ onto:

\begin{equation}
\label{riemann}
\widetilde{\Pi}_{\nu_{G},k\nu_{T}}\left(x,y\right)=\frac{d_{\nu_{G}}}{(2\pi)^{d_{T}}}\int_{G}\int_{T}\chi_{\nu_{G}}(g^{-1})t^{-k\nu_{T}}\Pi\left(\tilde{\mu}^{G}_{g^{-1}}\circ\tilde{\mu}^{T}_{t^{-1}}\left(x\right),y\right)dtdg,
\end{equation}
\noindent
where $d_{\nu_{G}}=\dimensione{(V_{\nu_{G}})}$ and $dg$, $dt$ are the associated measure for $G$ and 
$T$ such that $\int_{G}dg=1$ and $\int_{T}dt=1$. We start considering the diagonal case, so we have:

\begin{equation}
\label{scented}
\begin{multlined}[t][12.5cm]
\widetilde{\Pi}_{\nu_{G},k\nu_{T}}\left(x,x\right)=d_{\nu_{G}}\int_{G}\int_{T}\chi_{\nu_{G}}(g^{-1})t^{-k}\Pi\left(\tilde{\mu}^{G}_{g^{-1}}\circ\tilde{\mu}^{T}_{t^{-1}}\left(x\right),x\right)dt dg\\
=\frac{d_{\nu_{G}}}{2\pi}\int_{G}\int_{(-\pi,+\pi)^{d_{T}}}\chi_{\nu_{G}}(g^{-1})e^{-ik\nu_{T}\cdot\vartheta}\Pi\left(\tilde{\mu}^{G}_{g^{-1}}\circ\tilde{\mu}^{T}_{-\vartheta}\left(x\right),x\right)d\vartheta dg,
\end{multlined}
\end{equation}
\noindent
where $\vartheta\in(-\pi,\pi)^{d_{T}}$. For the moment suppose $x\in X$ generic and fixed, and denote $F_{x}\subseteq G\times T$ the stabilizer of $x$. 
For $\varepsilon>0$ we set 

$$ A=\{(g,t)\in G\times T: \distanza_{G\times T}{((g,t),F_{x})}<2\varepsilon\} $$
\noindent
and 

$$ B=\{(g,t)\in G\times T: \distanza_{G\times T}{((g,t),F_{x})}>\varepsilon\} $$
\noindent
so we have $G\times T=A\cup B$ and we can consider a partition of the unity $\gamma_{1}+\gamma_{2}=1$ associated to the covering $\{A,B\}$. 
We observe that the function:

\begin{equation}
\label{tea}
(g,t)\mapsto \gamma_{2}(g,t)\chi_{\nu_{G}}(g^{-1})\Pi\left(\tilde{\mu}^{G}_{g^{-1}}\circ\tilde{\mu}^{T}_{-\vartheta}\left(x\right),x\right)
\end{equation}
\noindent
is $\mathcal{C}^{\infty}$ because the singular support of $\Pi$ is included in the diagonal of $X\times X$. 

Then

\begin{equation}
\label{tea}
t\mapsto \int_{G}\gamma_{2}(g,t)\chi_{\nu_{G}}(g^{-1})\Pi\left(\tilde{\mu}^{G}_{g^{-1}}\circ\tilde{\mu}^{T}_{-\vartheta}\left(x\right),x\right)dg
\end{equation}
\noindent
is infinitely smooth and the Fourier transform is rapidly decreasing. Thus the contribution coming from $B$ is rapidly decreasing and we can multiply the integrand by $\gamma_{1}$. So we can only consider:

\begin{equation}
\label{store}
\begin{multlined}[t][12.5cm]
\widetilde{\Pi}_{\nu_{G},k\nu_{T}}\left(x,x\right)\\\sim \frac{d_{\nu_{G}}}{2\pi}\int_{G}\int_{(-\pi,\pi)^{d_{T}}}\gamma_{1}(g,\vartheta)\chi_{\nu_{G}}(g^{-1})e^{-ik\nu_{T}\cdot\vartheta}\Pi\left(\tilde{\mu}^{G}_{g^{-1}}\circ\tilde{\mu}^{T}_{-\vartheta}\left(x\right),x\right)d\vartheta dg.
\end{multlined}
\end{equation}

Now if $\gamma_{2}(g,\vartheta)\not=0$ then $\tilde{\mu}^{G}_{g^{-1}}\circ\tilde{\mu}^{T}_{-\vartheta}(x),x$ are near and we can represent $\Pi$ as 
Fourier integral operator as in $\cite{secondo}$:

\begin{equation}
\label{scotch}
\Pi(y,y')=\int_{0}^{+\infty}e^{it\psi(y,y')}s(y,y',t)dt,
\end{equation}
\noindent
where $\Im{(\psi)}\geq 0$ and $s$ is a semiclassical symbol admitting an asymptotic expansion $s(y,y',t)\sim\sum_{j=0}^{+\infty}t^{n-j}s_{j}(y,y')$.
Inserting $(\ref{scotch})$ in $(\ref{store})$ we obtain: 

\begin{equation}
\label{sharktale}
\begin{multlined}[t][12.5cm]
\widetilde{\Pi}_{\nu_{G},k\nu_{T}}\left(x,x\right)\\\sim \frac{d_{\nu_{G}}}{2\pi}\int_{G}\int_{(-\pi,+\pi)^{d_{T}}}\int_{0}^{+\infty}\gamma_{1}(g,\vartheta)\chi_{\nu_{G}}(g^{-1})e^{i[t\psi(\tilde{\mu}^{G}_{g^{-1}}\circ\tilde{\mu}^{T}_{-\vartheta}\left(x\right),x)-k\nu_{T}\cdot\vartheta]}\cdot\\
\cdot s\left(\tilde{\mu}^{G}_{g^{-1}}\circ\tilde{\mu}^{T}_{-\vartheta}\left(x\right),x,t\right) dt d\vartheta dg,
\end{multlined}
\end{equation}
\noindent
and performing the change of variables $t\rightarrow kt$, we get:

\begin{equation}
\label{sharktale2}
\begin{multlined}[t][12.5cm]
\widetilde{\Pi}_{\nu_{G},k\nu_{T}}\left(x,x\right)\sim \frac{d_{\nu_{G}}}{2\pi}k\cdot\int_{G}\int_{(-\pi,+\pi)^{d_{T}}}\int_{0}^{+\infty}\gamma_{1}(g,\vartheta)\chi_{\nu_{G}}(g^{-1})e^{ik\Psi(t,g,\vartheta,x)}\cdot\\
\cdot s\left(\tilde{\mu}^{G}_{g^{-1}}\circ\tilde{\mu}^{T}_{-\vartheta}\left(x\right),x,kt\right) dt d\vartheta dg,
\end{multlined}
\end{equation}
\noindent 
where we have set $\Psi(t,g,\vartheta,x)=t\psi(\tilde{\mu}^{G}_{g^{-1}}\circ\tilde{\mu}^{T}_{-\vartheta}\left(x\right),x)-\nu_{T}\cdot\vartheta$. We shall now use integration by parts in $\vartheta$ to prove that only a rapidly decreasing contribution to the asymptotic is lost, if the integrand in $(\ref{sharktale2})$
is multiplied  by a suitable cut-off function. In local coordinates we have $\tilde{\mu}^{G}_{g^{-1}}\circ\tilde{\mu}^{T}_{-\vartheta}\left(x\right)= x+O(\varepsilon)$
with $\varepsilon>0$ very small, because $\left(g,e^{i\vartheta}\right)\in U$ with $U$ a small neighborhood of $F_{x}$. Thus we have that:

$$d_{(\tilde{\mu}^{G}_{g^{-1}}\circ\tilde{\mu}^{T}_{-\vartheta}\left(x\right),x)}\psi=d_{(x,x)}\psi+O(\varepsilon)=\left(\alpha_{x},-\alpha_{x}\right)+O(\varepsilon),$$
\noindent
with $\partial_{\vartheta}\Psi= t\Phi_{T}(m)-\nu_{T}+O(\varepsilon)$. Therefore, since $\Phi_{T}(m)\not=\mathbf{0}$ and $\nu_{T}\not=\mathbf{0}$ we have for $t\gg 0$ that 

$$\|\partial_{\vartheta}\Psi\| \geq Ct,$$
\noindent
for some $C>0$. In a similar way for $0<t\ll 1$ we have 

$$\|\partial_{\vartheta}\Psi\|\geq C_{1}>0,$$
\noindent
for some $C_{1}>0$. Therefore by integration by parts in $d\vartheta$, we have  that the asymptotics for $k\rightarrow +\infty$ is unchanged. We multiply the integrand by $\rho(t)$, where $\rho\in\mathcal{C}^{\infty}_{0}\left(\frac{1}{2D},2D\right)$ and $\rho \equiv 1$ on $\left(\frac{1}{D},D\right)$, so that the integral in $dt$ is now compactly supported. 
 We shall now use integration by parts in $dt$ to show that only a rapidly decreasing contribution is lost, if the integration in $(g,\vartheta)$ is restricted to a
tubular neighborhood of $F_{x}$ of radius $O(k^{\delta-\frac{1}{2}})$. We have that $\partial_{t}\Psi=\psi(\tilde{\mu}^{G}_{g^{-1}}\circ\tilde{\mu}^{T}_{-\vartheta}\left(x\right),x)$. If $\distanza{\left(\left(g,e^{i\vartheta}\right),F_{x}\right)}\geq Ck^{\delta-\frac{1}{2}}$, then
$$\distanza{\left(\tilde{\mu}^{G}_{g^{-1}}\circ\tilde{\mu}^{T}_{-\vartheta}\left(x\right),x\right)}\geq Ck^{\delta-\frac{1}{2}}$$  
\noindent
and so:

\begin{equation}
 \label{intttt}
 \left|\psi(\tilde{\mu}^{G}_{g^{-1}}\circ\tilde{\mu}^{T}_{-\vartheta}\left(x\right),x)\right|\geq \Im{\psi(\tilde{\mu}^{G}_{g^{-1}}\circ\tilde{\mu}^{T}_{-\vartheta}\left(x\right),x)}\geq C_{2}k^{2\delta-1}
\end{equation}
\noindent
(see Corollary 2.3 of $\cite{secondo}$). Introducing the operator

$$L_{t}=\left[\psi(\tilde{\mu}^{G}_{g^{-1}}\circ\tilde{\mu}^{T}_{-\vartheta}\left(x\right),x)\right]^{-1}\partial_{t},$$
\noindent
we have that:

$$e^{ik\Psi}=-\frac{i}{k}L_{t}\left(e^{ik\Psi}\right).$$

We can now mimick the standard proof of the Stationary Phase Lemma: iteratively integrating by parts, we obtain at each step in view of $(\ref{intttt})$ a factor of order $O\left(k^{-2\delta}\right)$, and then after $N$ steps a factor of order $O\left(k^{-2N\delta}\right)$.
This proves that the contribution to the asymptotics  coming from the locus where $\distanza{\left(\left(g,e^{i\vartheta}\right),F_{x}\right)}\geq Ck^{\delta-\frac{1}{2}}$ is rapidly decreasing. We can now prove that $(\ref{sharktale2})$ is rapidly decreasing in $k$ for $\distanza_{X}{(x,X_{0,\nu_{T}})}\geq Ck^{\delta-\frac{1}{2}}$.
Now we consider a bump function $\rho_{1}:P\rightarrow \mathbb{R}$ supported in a small neighborhood of $F_{x}$ and $\equiv 1$ near to $F_{x}$. The function $\rho_{1}$ is defined as $\rho_{1}=\rho_{1}(f,\xi)$ with $f\in F_{x}$ and $\xi$ the normal coordinate to $F_{x}$. We can multiply the integrand of $(\ref{sharktale2})$ by $\rho_{1}\left(f,k^{\frac{1}{2}-\delta}\xi\right)$ losing only an $O(k^{-\infty})$. Then if $\rho_{1}(g,t)\not=0$ we have $\tilde{\mu}^{G}_{g^{-1}}\circ\tilde{\mu}^{T}_{-\vartheta}\left(x\right)=x + O\left(k^{\delta-\frac{1}{2}}\right)$. Therefore: 

$$d_{\left(\tilde{\mu}^{G}_{g^{-1}}\circ\tilde{\mu}^{T}_{-\vartheta}(x),x\right)}\psi
=d_{(x,x)}\psi+O\left(k^{\delta-\frac{1}{2}}\right)=\left(\alpha_{x},-\alpha_{x}\right)+O(k^{\delta-\frac{1}{2}}),$$

\noindent
and $\partial_{(\varsigma,\vartheta)}\psi=t\Phi_{P}(m)-\nu_{T}+ O\left(k^{\delta-\frac{1}{2}}\right)$. Here $(\varsigma,\vartheta)$ are local coordinates on $P$ induced by the exponential map $\exp_{P}$. Then if $\distanza_{X}{\left(x,X_{0,\nu_{T}}\right)}\geq C_{3}k^{\delta-\frac{1}{2}}$ we have that:

$$\|\partial_{(\varsigma,\vartheta)}\psi\|\geq C_{4}k^{\delta-\frac{1}{2}}.$$

Thus we find a differential operator $L_{\varsigma,\vartheta}$ with $\left|L_{\varsigma,\vartheta}\right|\geq C_{5}k^{\delta-\frac{1}{2}}$ where $\distanza_{X}{\left(x,X_{0,\nu_{T}}\right)}\geq O\left(k^{\delta-\frac{1}{2}}\right)$ such that $L_{\varsigma,\vartheta}\left(e^{ik\Psi}\right)=ike^{ik\Psi}$. Iterating the integration by parts, in view of the scaling factor we have at each step a factor $O\left(k^{-2\delta}\right)$. This proves that $\widetilde{\Pi}_{\nu_{G},k\nu_{T}}(x,x)=O(k^{-\infty})$ for $\distanza_{X}{\left(x,X_{0,\nu_{T}}\right)}\geq Ck^{\delta-\frac{1}{2}}$. Let us consider $(x,y)\in X\times X$ with 

$$\max{\{\distanza_{X}(x,X_{0,\nu_{T}}),\distanza_{X}(y,X_{0,\nu_{T}})\}}\geq Ck^{\delta-\frac{1}{2}}$$ 
\noindent
for every $\delta$ fixed and using the Cauchy-Schwarz inequality we have:

\begin{equation}
\label{antonio}
\left|\widetilde{\Pi}_{\nu_{G},k\nu_{T}}\left(x,y\right)\right|\leq \sqrt{\widetilde{\Pi}_{\nu_{G},k\nu_{T}}\left(x,x\right)}\cdot\sqrt{\widetilde{\Pi}_{\nu_{G},k\nu_{T}}\left(y,y\right)},
\end{equation}
\noindent
so $\widetilde{\Pi}_{\nu_{G},k\nu_{T}}\left(x,y\right)=O(k^{-\infty})$. This complete the proof of $1)$.

Let us now consider the proof of $2)$. Now setting $x_{jk}=x+\frac{(\theta_{j},\vv_{j})}{\sqrt{k}}$ for $j=1,2$, using FIO representation as before and changing variables $t\rightarrow kt$, we get:

\begin{equation}
\label{newton}
\begin{multlined}[t][12.5cm]
\widetilde{\Pi}_{\nu_{G},k\nu_{T}}\left(x_{1k},x_{2k}\right)=\frac{kd_{\nu_{G}}}{(2\pi)^{d_{T}}}\int_{W}\chi_{\nu_{G}}(g^{-1})e^{ik\Psi^{(1)}(t,\vartheta,x)}\cdot\\
\cdot s\left(\tilde{\mu}^{G}_{g^{-1}}\circ\tilde{\mu}^{T}_{-\vartheta}\left(x_{1k}\right),x_{2k},kt\right)dV_{W}(w),
\end{multlined}
\end{equation}
\noindent
where 

\begin{equation}
W=G\times (-\pi,\pi)^{d_{T}}\times (0,+\infty) \ \ \ \ dV_{W}(w)=dg d\vartheta dt,
\end{equation}
\noindent
and

\begin{equation}
\label{jacobi}
\Psi^{(1)}(t,\vartheta,x)=t\psi\left(\tilde{\mu}^{G}_{g^{-1}}\circ\tilde{\mu}^{T}_{-\vartheta}\left(x_{1k}\right),x_{2k}\right)-\nu_{T}\cdot \vartheta.
\end{equation}

Here $t=(t_{1}, \cdots, t_{d_{T}})=\left(e^{i\vartheta_{1}},\cdots , e^{i\vartheta_{d_{T}}}\right)=e^{i\vartheta}$. Let $F_{m}\subseteq P$, $F_{m}=\{p_{j}\}=\{(g_{j},t_{j})\}$ the finite stabilizer of $x\in X_{0,\nu_{T}}$. We introduce a bump function $\rho=\sum_{j=1}^{N_{x}}\rho_{j}$ with support of $\rho_{j}$ in a neighborhood of $p_{j}=(g_{j},t_{j})$. As consequence we have:

\begin{equation}
\label{eulero}
\widetilde{\Pi}_{\nu_{G},k\nu_{T}}\left(x_{1k},x_{2k}\right)\sim \sum_{j}\Pi_{\nu_{G},k\nu_{T}}\left(x_{1k},x_{2k}\right)^{(j)},
\end{equation}
\noindent
where each addend of $(\ref{eulero})$ is given by $(\ref{newton})$ multiplied by $\rho_{j}$. 
In the support of each $p_{j}$ we write $g=g_{j}\exp_{G}{\left(\frac{\gamma}{\sqrt{k}}\right)}$ and $t=t_{j}e^{\frac{i\vartheta}{\sqrt{k}}}$, where with $\exp_{G}$ we denote the exponential map from $\mathfrak{g}\rightarrow G$ and $\gamma,\vartheta$ are coordinates respectively on $\mathfrak{g}\cong \mathbb{R}^{d_{G}}$,$\mathfrak{t}\cong \mathbb{R}^{d_{T}}$ associated with the respective orthonormal basis. Omitting $\rho_{j}$ in the integrand we have:

\begin{equation}
\label{leibnitz}
\begin{multlined}[t][12.5cm]
\widetilde{\Pi}_{\nu_{G},k\nu_{T}}\left(x_{1k},x_{2k}\right)\sim\frac{k^{1-\frac{d_{P}}{2}}d_{\nu_{G}}}{(2\pi)^{d_{T}}}\int_{W'}\chi_{\nu_{G}}\left(g^{-1}_{j}\exp_{G}{\left(-\frac{\gamma}{\sqrt{k}}\right)}\right)t_{j}^{-1}e^{ik\Psi^{(2)}(t,\vartheta,x)}\cdot\\ \\
\cdot s\left(\tilde{\mu}^{G}_{-\frac{\gamma}{\sqrt{k}}}\circ\tilde{\mu}^{T}_{-\frac{\vartheta}{\sqrt{k}}}\left(x_{1k}\right),x_{2k},kt\right)dV_{W'}(w),
\end{multlined}
\end{equation}
\noindent
where 

\begin{equation}
W'=\mathbb{R}^{d_{G}}\times \mathbb{R}^{d_{T}}\times (0,+\infty) \ \ \ \ dV_{W'}(w)=d\gamma d\vartheta dt
\end{equation}
\noindent
and

\begin{equation}
\label{cantor}
\Psi^{(2)}(t,\vartheta,x)=t\psi\left(\tilde{\mu}^{G}_{-\frac{\gamma}{\sqrt{k}}}\circ\tilde{\mu}^{T}_{-\frac{\vartheta}{\sqrt{k}}}\circ\tilde{\mu}^{G}_{g^{-1}_{j}}\circ\tilde{\mu}^{T}_{t^{-1}_{j}}\left(x_{1k}\right),x_{2k}\right)-\nu_{T}\cdot \frac{\vartheta}{\sqrt{k}}.
\end{equation}

We write $x_{1k}^{j}=\tilde{\mu}^{G}_{g^{-1}_{j}}\circ\tilde{\mu}^{T}_{t^{-1}_{j}}\left(x_{1k}\right)=x+\frac{1}{\sqrt{k}}(\theta_{1},\vv_{1}^{j})$, for a particular choice of $\sigma$ adapted section in the definition of HLC. We assume that the orthonormal basis of $\mathfrak{t}$ is taken as $(w_{1},\cdots,w_{d_{T}})$ with $(w_{1},\cdots,w_{d_{T}-1})$ an orthonormal basis for $\kernel{(\Phi_{T}(m))}$ and $\langle \Phi_{T}(m),w_{d_{T}}\rangle=\|\Phi_{T}(m)\|$.
So we have that, if $(v_{1},\cdots ,v_{d_{G}})$ is the orthonormal basis for $\mathfrak{g}$, an orthonormal basis for $\mathfrak{p}=\lie{(P)}$ is of the form:

$$ (v_{1},\cdots,v_{d_{G}},w_{1},\cdots,w_{d_{T}}=\eta). $$

We call $(a_{1},\cdots, a_{d_{P}}=b)$ the corresponding linear coordinates on $\mathfrak{p}$ such that $a=(a_{1},\cdots, a_{d_{P}}=b)\in \mathbb{R}^{d_{P}-1}\cong \kernel{(\Phi_{P}(m))}$ and $a_{M}(m)\in \mathbb{R}^{2d_{M}}\cong T_{m}M$ is his injective valutation. Considering $(\ref{cantor})$, we write $(\gamma,\vartheta)\in \mathfrak{p}\cong \mathbb{R}^{d_{P}}\cong \mathbb{R}^{d_{G}}\times \mathbb{R}^{d_{T}}$ as

\begin{equation}
\label{zermelo}
(\gamma,\vartheta)=(\gamma,\vartheta')+\vartheta\eta=a+b\eta,
\end{equation}
\noindent
remember that $\nu_{T}=\lambda\cdot \Phi_{T}(m) \Rightarrow (\mathbf{0},\nu_{T})=\lambda\Phi_{P}(m)$, then we have:

$$\nu_{T}\cdot \vartheta=(\mathbf{0},\nu_{T})\cdot (a+b\eta)=b\lambda\|\Phi_{T}(m)\|.$$

Here $\lambda=\lambda_{\nu_{T}}$ is such that $\nu_{T}=\lambda_{\nu_{T}}\cdot\Phi_{T}(m)$. Thus we have:

\begin{equation}
\label{hilbert}
\Psi^{(2)}(t,b,x)=t\psi\left(\tilde{\mu}^{P}_{-\frac{a+b\eta}{\sqrt{k}}}\left(x_{1k}^{j}\right),x_{2k}\right)-\lambda_{\nu_{T}}\frac{\|\Phi_{T}(m)\|b}{\sqrt{k}}.
\end{equation}

Now let $p=a+b\eta$ we have that $a\in\kernel{(\Phi_{P}(m))}$ and $\eta\in\kernel{(\Phi_{P}(m))}^{\perp}$, with $\|\eta\|=1$, $\langle\Phi_{P}(m),\eta\rangle=\langle\Phi_{T}(m),\eta\rangle=\|\Phi_{T}(m)\|$.
We get:

\begin{equation}
\label{godel}
\begin{multlined}[t][12.5cm]
\tilde{\mu}^{P}_{-\frac{p}{\sqrt{k}}}\left(x_{1k}^{j}\right)=\tilde{\mu}^{P}_{-\frac{p}{\sqrt{k}}}\left(x+\frac{1}{\sqrt{k}}(\theta_{1},\vv_{1}^{j})\right)=\\ \\
=x+\left(\frac{\theta_{1}+\langle\Phi_{T}(m),p\rangle}{\sqrt{k}}+\frac{1}{k}\omega_{m}(a_{M}(m)+b\eta_{M}(m),\vv_{1}^{j})+ \widetilde{B}_{3}\left(\frac{a_{M}(m)}{\sqrt{k}},\frac{\vv}{\sqrt{k}},\frac{b}{\sqrt{k}}\right),\right.\\
\left.,\frac{\theta_{1}}{\sqrt{k}}(\vv_{1}^{j}-a_{M}(m)-b\eta_{M})+\widetilde{B}_{2}\left(\frac{a_{M}(m)}{\sqrt{k}},\frac{\vv}{\sqrt{k}},\frac{b}{\sqrt{k}}\right)\right)
\\
\\=x+(\mathcal{A}_{j,k}(\vartheta,v_{1}^{j},v_{2}),\mathcal{B}_{j,k}(\vartheta,v_{1}^{j},v_{2})),
\end{multlined}
\end{equation}
\noindent
with $\widetilde{B}_{2},\widetilde{B}_{3}$ that vanish at the origin to third and second order. 
We obtain that:\\
\\
\noindent
$t\psi\left(\tilde{\mu}^{P}_{-\frac{p}{\sqrt{k}}}\left(x_{1k}^{j}\right),x_{2k}\right)-\frac{\vartheta}{\sqrt{k}}$
\begin{equation}
\label{einstein}
\begin{multlined}[t][12.5cm]
 =it\left[1-e^{i\left(\mathcal{A}_{j,k}(\vartheta,v_{1}^{j},v_{2})-\frac{\theta_{2}}{\sqrt{k}}\right)}\right]\\ -\frac{it}{k}\psi_{2}\left(\vv_{1}^{j}-a_{M}(m)-b\eta_{M}(m),\vv_{2}\right)e^{i\left(\mathcal{A}_{j,k}(\vartheta,v_{1}^{j},v_{2})-\frac{\theta_{2}}{\sqrt{k}}\right)}-\frac{\vartheta}{\sqrt{k}}+\\
 +it R_{3}^{\psi}\left(\frac{1}{\sqrt{k}}\left(\vv_{1}^{j}-\vartheta\xi_{M}(m)-\varsigma_{M}(m)\right),\frac{\vv_{2}}{\sqrt{k}}\right)e^{i\left(\mathcal{A}_{j,k}(\vartheta,v_{1}^{j},v_{2})-\frac{\theta_{2}}{\sqrt{k}}\right)},
\end{multlined}
\end{equation}
\noindent
where $R_{3}^{\psi}$ vanishes to third order at the origin and

$$ \psi_{2}(r,s)=-i\omega_{m}(r,s)-\frac{1}{2}\|r-s\|^2  \ \ \ \ (r,s\in\mathbb{C}^{n}).$$

Now

\begin{equation*}
\begin{multlined}[t][12.5cm]
i\left(\mathcal{A}_{j,k}(\vartheta,v_{1}^{j},v_{2})-\frac{\theta_{2}}{\sqrt{k}}\right)=\\
=\frac{i}{\sqrt{k}}(\theta_{1}-\theta_{2}+b\|\Phi_{T}\|)+\frac{i}{k}\omega_{m}(a_{M}+b\eta_{M},\vv_{1}^{j})+ \widetilde{B}_{3}'\left(\frac{a_{M}(m)}{\sqrt{k}},\frac{\vv}{\sqrt{k}},\frac{b}{\sqrt{k}}\right), 
\end{multlined}
\end{equation*}
\noindent
then 

\begin{equation*}
\begin{multlined}[t][12.5cm]
1-e^{i\left(\mathcal{A}_{j,k}(\vartheta,v_{1}^{j},v_{2})-\frac{\theta_{2}}{\sqrt{k}}\right)}=\\
=1-\left\{ 1+ \frac{i}{\sqrt{k}}(\theta_{1}-\theta_{2}+b\|\Phi_{T}\|)+ \frac{i}{k}\omega_{m}\left(a_{M}(m)+b\eta_{M}(m),\vv_{1}^{j}\right)-\right.\\
\left.-\frac{1}{2k}(\theta_{1}-\theta_{2}+b\|\Phi_{T}\|)^2+\widetilde{B}_{3}''\left(\frac{a_{M}(m)}{\sqrt{k}},\frac{\vv}{\sqrt{k}},\frac{b}{\sqrt{k}}\right)\right\}\\
=- \frac{i}{\sqrt{k}}(\theta_{1}-\theta_{2}+b\|\Phi_{T}\|)-\frac{i}{k}\omega_{m}\left(a_{M}(m)+b\eta_{M}(m),\vv_{1}^{j}\right)+\\
+\frac{1}{2k}(\theta_{1}-\theta_{2}+b\|\Phi_{T}\|)^2+\widetilde{B}_{3}'''\left(\frac{a_{M}(m)}{\sqrt{k}},\frac{\vv}{\sqrt{k}},\frac{b}{\sqrt{k}}\right)
\end{multlined}
\end{equation*}
\noindent
and 

\begin{equation*}
\begin{multlined}[t][12.5cm]
it\left[1-e^{i\left(\mathcal{A}_{j,k}(\vartheta,v_{1}^{j},v_{2})-\frac{\theta_{2}}{\sqrt{k}}\right)}\right]=\\
=\frac{t}{\sqrt{k}}(\theta_{1}-\theta_{2}+b\|\Phi_{T}\|)+\frac{1}{k}t\omega_{m}\left(a_{M}(m)+b\eta_{M}(m),\vv_{1}^{j}\right)+\\
+\frac{i}{2k}t(\theta_{1}-\theta_{2}+b\|\Phi_{T}\|)^2+it\widetilde{B}_{3}'''\left(\frac{a_{M}(m)}{\sqrt{k}},\frac{\vv}{\sqrt{k}},\frac{b}{\sqrt{k}}\right).
\end{multlined}
\end{equation*}

Thus as conseguence:

\begin{equation*}
\label{weyl}
\begin{multlined}[t][12.5cm]
\Psi^{(2)}(t,b,x)=\frac{t}{\sqrt{k}}(\theta_{1}-\theta_{2}+b\|\Phi_{T}\|)+\\
+\frac{1}{k}t\omega_{m}\left(a_{M}(m)+b\eta_{M}(m),\vv_{1}^{(j)}\right)+\frac{i}{2k}t(\theta_{1}-\theta_{2}+b\|\Phi_{T}\|)^2-\lambda_{\nu_{T}}\frac{\|\Phi_{T}(m)\|b}{\sqrt{k}}-\\
-\frac{it}{k}\psi_{2}\left(\vv_{1}^{j}-a_{M}(m)-b\eta_{M}(m),\vv_{2}\right)+it\widetilde{B}_{3}'''\left(\frac{a_{M}(m)}{\sqrt{k}},\frac{\vv}{\sqrt{k}},\frac{b}{\sqrt{k}}\right)\\
\end{multlined}
\end{equation*}
\begin{equation*}
\begin{multlined}[t][12.5cm]
=\frac{1}{\sqrt{k}}\left(t(\theta_{1}-\theta_{2}+b\|\Phi_{T}\|)-\lambda_{\nu_{T}}\|\Phi_{T}(m)\|b\right)+\\
+\frac{1}{k}\left[t\omega_{m}\left(a_{M}(m)+b\eta_{M}(m),\vv_{1}^{(j)}\right)+\frac{i}{2}t(\theta_{1}-\theta_{2}+b\|\Phi_{T}\|)^2-\right.\\
\left.-it\psi_{2}\left(\vv_{1}^{j}-a_{M}(m)-b\eta_{M}(m),\vv_{2}\right) \right]+it\widetilde{B}_{3}'''\left(\frac{a_{M}(m)}{\sqrt{k}},\frac{\vv}{\sqrt{k}},\frac{b}{\sqrt{k}}\right)
\end{multlined}
\end{equation*}
\noindent
and \\
\\
\noindent
$ik\Psi^{(2)}(t,b,x)$
\begin{equation*}
\label{weyl}
\begin{multlined}[t][12.5cm]
=i\sqrt{k}\left(t(\theta_{1}-\theta_{2}+b\|\Phi_{T}\|)-\lambda_{\nu_{T}}\|\Phi_{T}(m)\|b\right)+\\
+\left[it\omega_{m}\left(a_{M}(m)+b\eta_{M}(m),\vv_{1}^{(j)}\right)-\frac{1}{2}t(\theta_{1}-\theta_{2}+b\|\Phi_{T}\|)^2+\right.\\
\left.+t\psi_{2}\left(\vv_{1}^{j}-a_{M}(m)-b\eta_{M}(m),\vv_{2}\right) \right]-kt\widetilde{B}_{3}'''\left(\frac{a_{M}(m)}{\sqrt{k}},\frac{\vv}{\sqrt{k}},\frac{b}{\sqrt{k}}\right).
\end{multlined}
\end{equation*}

Continuing calculations we can rewrite the $j$-term as:

\begin{equation}
\label{weil}
\begin{multlined}[t][12.5cm]
\widetilde{\Pi}_{\nu_{G},k\nu_{T}}(x_{1k},x_{2k})^{(j)} \\
\sim k^{1-\frac{d_{P}}{2}}\frac{d_{\nu_{G}}}{(2\pi)^{d_{T}}}\cdot\int_{\mathbb{R}^{d_{P}-1}}da\cdot\left[\int_{0}^{+\infty}dt\int_{-\infty}^{+\infty}db e^{i\sqrt{k}\Upsilon(t,b)}e^{A(m,\theta,\vv,v,a)}B(j)\right],
\end{multlined}
\end{equation}
\noindent
where

\begin{equation}
\label{dedekind}
\Upsilon(t,b)=t\left(b\|\Phi_{T}(m)\|+\theta_{1}-\theta_{2}\right)-\lambda\|\Phi_{T}\|b,
\end{equation}

\begin{equation}
\begin{multlined}[t][12.5cm]
A(m,\theta,\vv,v,a)\\
=-\frac{t}{2}(b\|\Phi_{T}(m)\|+\theta_{1}-\theta_{2})^2+it\omega_{m}(a_{M}(m)+b\eta_{M}(m),\vv_{1}^{j})+\\
+t\psi_{2}(\vv_{1}^{j}-a_{M}(m)-b\eta_{M}(m),\vv_{2})e^{i\mathcal{A}_{j,k}(\vartheta,v_{1}^{j},v_{2})} 
\end{multlined}
\end{equation}
\noindent
and

\begin{equation}
\label{degiorgi}
B(j)=\chi_{\nu_{G}}(g_{j}^{-1}) e^{-ik\vartheta_{j}\nu_{T}}.
\end{equation}

The internal integral in $(\ref{weil})$ is oscillatory in $\sqrt{k}$ with phase $\Upsilon$. The phase has critical points $(t_{0},b_{0})=\left(\lambda_{\nu_{T}},\frac{\theta_{2}-\theta_{1}}{\|\Phi_{T}(m)\|}\right)$. The Hessian is:

$$H(\Upsilon)(P_{0})=\left(\begin{array}{cc} 0 & \|\Phi_{T}(m)\| \\ \|\Phi_{T}(m)\| & 0 \end{array}\right)$$ 
\noindent
and 

$$\sqrt{\left|\det{\frac{\sqrt{k}H}{2\pi}}\right|}=\frac{\sqrt{k}}{2\pi}\|\Phi_{T}(m)\|.$$

Using the Stationary Phase Lemma $\cite{settimo}$ we have:

\begin{equation}
\label{pitagora}
\begin{multlined}[t][12.5cm]
\widetilde{\Pi}_{\nu_{G},k\nu_{T}}(x_{1k},x_{2k})^{(j)}\sim \frac{d_{\nu_{G}}}{(2\pi)^{d_{T}-1}}B(j)k^{\frac{1}{2}-\frac{d_{P}}{2}}\cdot\int_{\mathbb{R}^{d_{P}-1}}e^{-i\sqrt{k}(\theta_{2}-\theta_{1})\lambda_{\nu}}\cdot\\
\cdot e^{A(a,t_{0},b_{0})}\cdot\frac{1}{\|\Phi_{T}(m)\|}\left(\frac{k\lambda_{\nu_{T}}}{\pi}\right)^{d_{M}}da,
\end{multlined}
\end{equation}
\noindent
where 

\begin{equation*}
\begin{multlined}[t][12.5cm]
A(a,t_{0},b_{0})=\lambda_{\nu_{T}}\left[i\omega_{m}\left(a_{M}(m)+b_{0}\eta_{M},\vv_{1}^{j}\right)+\right.\\
\left.+\psi_{2}\left(\vv_{1}^{j}-a_{M}(m)-b_{0}\eta_{M}(m),\vv_{2}\right)\right],
\end{multlined}
\end{equation*}
\noindent
and we have to evaluate the integral $\int_{\mathbb{R}^{d_{P}-1}}da e^{A(a,t_{0},b_{0})}$. In order to do this we define the following spaces:

$$V_{m}=\valutation_{m}{\left(\mathfrak{g}\oplus\kernel{(\Phi_{T}(m))}\right)}=\valutation_{m}{(\kernel{(\Phi_{P}(m))})} \subseteq T_{m}M,$$
$$N_{m}= J_{m}(V_{m}),$$
$$H_{m}=[V_{m}\oplus N_{m}]^{\perp},$$
\noindent
where $V_{m},H_{m}$ are complex subspaces of $T_{m}M$ and $N_{m}$ is the normal space to $M_{0,\nu_{T}}$. Decomposing $\vv\in T_{m}M$ as $\vv=\vv_{h}+\vv_{v}+\vv_{t}$, with respectively $\vv_{h}\in H_{m}$, $\vv_{v}\in V_{m}$ and $\vv_{t}\in N_{m}$, we have that $a_{M}=a_{Mv}$, $\eta_{M}=\eta_{Mh}+\eta_{Mv}$ and\\
\\
\noindent
$i\omega_{m}\left(a_{M}(m)+b_{0}\eta_{M}(m),\vv_{1}^{j}\right)+i\psi_{2}\left(\vv_{1}^{j}-a_{M}(m)-b_{0}\eta_{M}(m),\vv_{2}\right)$
\begin{equation}
\label{fermat2}
\begin{multlined}[t][12.5cm]
=i\omega_{m}(a_{M}(m),\vv_{1t}^{j}+\vv_{2t})-\frac{1}{2}\|\vv_{1v}-a_{M}(m)-b_{0}\eta_{Mv}(m)-\vv_{2v}\|^{2}\\
+ib_{0}\omega_{m}(\eta_{Mh}(m),\vv_{1h}^{j}+\vv_{2h})+ib_{0}\omega_{m}(\eta_{Mv}(m),\vv_{1t}^{j}+\vv_{2t})-i\omega_{m}(\vv_{1h}^{j},\vv_{2h})+\\
-\frac{1}{2}\|\vv_{1h}^{j}-b_{0}\eta_{Mh}(m)-\vv_{2h}\|^{2}-i\omega_{m}(\vv_{1v}^{j},\vv_{2t})-\frac{1}{2}\|\vv_{1t}^{j}-\vv_{2t}\|^{2}-i\omega_{m}(\vv_{1t}^{j},\vv_{2v}).
\end{multlined}
\end{equation}

We define $r_{M}(m)\in \mathbb{R}^{d_{P}-1}$ translated by $a_{M}(m)$ such that: 

$$ r_{M}(m)=\vv_{1v}^{j}-a_{M}(m)-b_{0}\eta_{Mv}(m)-\vv_{2v}, $$
\noindent
so $ a_{M}(m)=\vv_{1v}^j-r_{M}(m)-b_{0}\eta_{Mv}(m)-\vv_{2v}$ and we get:

\begin{equation}
\label{gaylussac2}
\begin{multlined}[t][12.5cm]
\omega_{m}(a_{M}(m),\vv_{1t}^{j}+\vv_{2t})=\omega_{m}(\vv_{1v}^{j},\vv_{1t}^{j})+\omega_{m}(\vv_{1v}^{j},\vv_{2t})-\omega_{m}(\vv_{2v},\vv_{1t}^{j})-\omega_{m}(\vv_{2v},\vv_{2t})-\\
-\omega_{m}(r_{M}(m),\vv_{1t}^{j}+\vv_{2t})-\omega_{m}(b_{0}\eta_{Mv}(m),\vv_{1t}^{j}+\vv_{2t}).
\end{multlined}
\end{equation}

Putting $(\ref{gaylussac2})$ in $(\ref{fermat2})$ and deleting the opposite terms, we obtain:\\
\\
\noindent
$i\omega_{m}\left(a_{M}(m)+b_{0}\eta_{M}(m),\vv_{1}^{j}\right)+i\psi_{2}\left(\vv_{1}^{j}-a_{M}(m)-b_{0}\eta_{M}(m),\vv_{2}\right)$
\begin{equation*}
\begin{multlined}[t][12.5cm]
=i\left[\omega_{m}(\vv_{1v}^{j},\vv_{1t}^{j})-\omega_{m}(\vv_{2v},\vv_{2t})\right]-\frac{1}{2}\|\vv_{1t}^j-\vv_{2t}\|^{2} +\\
+ib_{0}\omega_{m}(\eta_{Mh}(m),\vv_{1h}^{j}+\vv_{2h})-i\omega_{m}(\vv_{1h}^{j},\vv_{2h})-\frac{1}{2}\left\|\vv_{1h}^{j}-b_{0}\eta_{Mh}(m)-\vv_{2h}\right\|^2-\\
-i\omega_{m}(r_{M}(m),\vv_{1t}^{j}+\vv_{2t})-\frac{1}{2}\|r_{M}(m)\|^2.
\end{multlined}
\end{equation*}

Let $C$ the matrix of $\valutation_{m}:\kernel{(\Phi_{T}(m))}\oplus\mathfrak{g}\rightarrow V_{m}$. Changing variable $r'=Cr$ we have that $dr=\det{C}^{-1}dr'$. We have to evaluate:

\begin{equation}
\label{nash}
\frac{1}{\det{C}}\int_{\mathbb{R}^{d_{P}-1}}e^{\lambda_{\nu_{T}}\left[-i\omega_{m}(r',\vv_{1t}^{j}+\vv_{2t})-\frac{1}{2}\|r'\|^2\right]}dr'.
\end{equation}

We make the substitution $s=\sqrt{\lambda_{\nu_{T}}}r'$, so we obtain:

\begin{equation}
\label{nash2}
=\frac{\lambda_{\nu_{T}}^{-\frac{1}{2}(d_{P}-1)}}{\det{C}}\int_{\mathbb{R}^{d_{P}-1}}e^{\left[-i\omega_{m}(s,\sqrt{\lambda_{\nu_{T}}}(\vv_{1t}^{j}+\vv_{2t}))-\frac{1}{2}\|s\|^2\right]}ds=\frac{(2\pi)^{\frac{d_{P}-1}{2}}}{\lambda^{\frac{1}{2}(d_{P}-1)}}e^{-\frac{1}{2}\lambda_{\nu_{T}}\|\vv_{1t}^{j}+\vv_{2t}\|^{2}}.
\end{equation}

Thus the exponential factor in the asymptotic expansion is $e^{H(v_{1}^{j},v_{2})}$ with 

\begin{equation}
\begin{multlined}[t][12.5cm]
H(v_{1}^{j},v_{2})=\lambda_{\nu_{T}}\left(-i\omega_{m}\left(\vv_{1h}^{j},\vv_{2h}\right)-\|\vv_{1t}\|^2-\|\vv_{2t}\|^2+ib_{0}\omega_{m}\left(\eta_{Mh}(m),\vv_{1h}^{j}+\vv_{2h}\right)\right.\\
\left.+i[\omega_{m}\left(\vv_{1v}^{j},\vv_{1t}^j\right)-\omega_{m}\left(\vv_{2v},\vv_{2t}\right)]-\frac{1}{2}\left\|\vv_{1h}^{j}-b_{0}\eta_{Mh}(m)-\vv_{2h}\right\|^2\right),
\end{multlined}
\end{equation}
\noindent
where $b_{0}=\frac{\theta_{2}-\theta_{1}}{\|\Phi_{T}(m)\|}$ and the principal term is of the form:

\begin{equation}
\label{colzanileonardo}
\begin{multlined}[t][12.5cm]
\frac{d_{\nu_{G}}2^{\frac{d_{G}}{2}}}{(\sqrt{2}\pi)^{d_{T}-1}}\left(\frac{k\lambda_{\nu_{T}}}{\pi}\right)^{d_{M}-\frac{d_{G}}{2}-\frac{(d_{T}-1)}{2}}\cdot\chi_{\nu_{G}}(g_{j}^{-1}) e^{-ik\vartheta_{j}\nu_{T}}\cdot\\
 \cdot\frac{e^{i\sqrt{k}(\theta_{2}-\theta_{1})\lambda_{\nu}}\cdot e^{H}}{\|\Phi_{T}(m)\|\mathcal{D}(m)}.
\end{multlined}
\end{equation}

This complete the proof of $2)$.

Let us now describe the necessary changes to the previous argument to prove $3)$. 
Instead of considering a neighborhood of $(x,x)$ we consider the asymptotic in a neighborhood of $(x,p_{0}\cdot x)$.
We set $y=p\cdot x$ and we assume given the local system of Heisenberg coordinates in a neighborhood of $x$ and $y$. We may also assume without loss that the Heisenberg coordinate system centered at $y$ is obtained from the one centered at $x$ by a $p$-translation, that is, 

$$y+(\theta,\vv)=p\cdot \big(x+(\theta,\vv)\big).$$

We must evaluate:

$$ \widetilde{\Pi}_{\nu_{G},k\nu_{T}}\left(x_{1k},y_{2k}\right),$$
\noindent
where $x_{1k}=x+\frac{u_{1}}{\sqrt{k}}$, $y_{2k}=y+\frac{u_{2}}{\sqrt{k}}$ and $u_{j}=\left(\theta_{j},\vv_{j}\right)$. Now with the preceding interpretation we have that $y_{2k}=p_{0}\left(x+\frac{u_{2}}{\sqrt{k}}\right)=p_{0}x_{2k}$. Proceeding as in the previous case we have:

\begin{equation}
\label{riemannnnn1}
\begin{multlined}[t][12.5cm]
\widetilde{\Pi}_{\nu_{G},k\nu_{T}}\left(x_{1k},y_{2k}\right)\\=\frac{d_{\nu_{G}}}{(2\pi)^{d_{T}}}\int_{G}\int_{T}\chi_{\nu_{G}}(g^{-1})t^{-k\nu_{T}} \cdot\Pi\left(\tilde{\mu}^{G}_{g^{-1}}\circ\tilde{\mu}^{T}_{t^{-1}}\left(x\right),p_{0}\cdot x_{2k}\right)dtdg,
\end{multlined}
\end{equation}
\noindent
and due to the unitary action:

\begin{equation}
\label{riemannnnn2}
\begin{multlined}[t][12.5cm]
\widetilde{\Pi}_{\nu_{G},k\nu_{T}}\left(x_{1k},y_{2k}\right)\\=\frac{d_{\nu_{G}}}{(2\pi)^{d_{T}}}\int_{G}\int_{T}\chi_{\nu_{G}}(g^{-1})t^{-k\nu_{T}}
\cdot\Pi\left(\tilde{\mu}^{P}_{p_{0}^{-1}}\circ\tilde{\mu}^{G}_{g^{-1}}\circ\tilde{\mu}^{T}_{t^{-1}}\left(x\right), x_{2k}\right)dtdg.
\end{multlined}
\end{equation}

Let $p_{0}=\left(g_{0},t_{0}\right)=\left(g_{0},e^{i\vartheta_{0}}\right)\in P$, then 

$$ \tilde{\mu}^{P}_{p_{0}^{-1}}\circ\tilde{\mu}^{P}_{(g,t)^{-1}}(x_{1k})=\tilde{\mu}_{g_{0}^{-1}g^{-1}}^{G}\circ\tilde{\mu}_{t_{0}^{-1}t^{-1}}^{T}(x_{1k}), $$
\noindent
and changing variables $g'=gg_{0}$, $t'=tt_{0}$, we have $g=g'g_{0}^{-1}$, $t=t't_{0}^{-1}$ and $\vartheta=\vartheta'-\vartheta_{0}$.

Then 
 
\begin{equation}
\label{riemannnnnstacippa}
\begin{multlined}[t][12.5cm]
\widetilde{\Pi}_{\nu_{G},k\nu_{T}}\left(x_{1k},y_{2k}\right)\\
=\frac{d_{\nu_{G}}}{(2\pi)^{d_{T}}}\int_{G}\int_{T}\overline{\chi_{\nu_{G}}(g'g_{0}^{-1})}e^{-ik\nu_{T}(\vartheta'-\vartheta_{0})}\Pi\left(\tilde{\mu}^{P}_{p_{0}^{-1}}\circ\tilde{\mu}^{G}_{g'^{-1}}\circ\tilde{\mu}^{T}_{t'^{-1}}\left(x\right), x_{2k}\right)dt'dg',
\end{multlined}
\end{equation}
\noindent
that will be the same as before with the difference that in the $j$-addendum we will make the substitution: $$\overline{\chi_{\nu_{G}}\left(g_{j}\right)}t_{j}^{-k\nu_{T}}=\overline{\chi_{\nu_{g},k\nu_{T}}(p_{j})}\mapsto \overline{\chi_{\nu_{G}}\left(g_{j}g_{0}^{-1}\right)}(t_{j}t_{0}^{-1})^{-k\nu_{T}}=\overline{\chi_{\nu_{G},k\nu_{T}}(p_{j}p_{0}^{-1})}.$$

This complete the proof of $3)$ and complete the proof of the main Theorem.
\hfill $\Box$

\section{Proof of Theorem $\bf{~\ref{teo:diagonal}}$}

$\Proof.$

On the diagonal we have:

\begin{equation}
\label{arnold}
\begin{multlined}[t][12.5cm]
\widetilde{\Pi}_{\nu_{G},k\nu_{T}}(x,x)\\
=d_{\nu_{G}}\int_{G}\int_{T}\chi_{\nu_{G}}(g^{-1})t^{-k\nu_{T}}\Pi\left(\tilde{\mu}_{g^{-1}}^{G}\circ\tilde{\mu}_{t^{-1}}^{T}(x),x\right)dt dg\\
=\frac{d_{\nu_{G}}}{(2\pi)^{d_{T}}}\int_{G}\int_{(-\pi,\pi)^{d_{T}}}\chi_{\nu_{G}}(g^{-1})e^{-ik\nu_{T}\cdot\vartheta}\Pi\left(\tilde{\mu}_{g^{-1}}^{G}\circ\tilde{\mu}_{-\vartheta}^{T}(x),x\right)d\vartheta dg,
\end{multlined}                                             
\end{equation}
\noindent
where $\vartheta=(\vartheta_{1},\cdots\cdots, \vartheta_{d_{T}})$, $\nu_{T}=(\nu_{1},\cdots,\nu_{d_{T}})\in\mathbb{Z}^{d_{T}}$ and $\nu_{T}\cdot\vartheta=\sum_{j=1}^{d_{T}}{\nu_{T}}_{j}\vartheta_{j}$. 
We consider $F_{x}=\left\{(g_{1},t_{1}),\cdots ,(g_{N_{x}},t_{N_{x}})\right\}$ the stabilizer, with $|F_{x}|=N_{x}$.
Let $\varepsilon >0$ and we consider the following open subsets of $P=G\times T$:

$$A=\{(g,t)\in G\times T : \distanza_{G\times T}{((g,t),F_{x})}<2\varepsilon\}$$

$$B=\{(g,t)\in G\times T : \distanza_{G\times T}{((g,t),F_{x})}>\varepsilon\}.$$

Then $P=A\cup B$, and we have choose a partition of unity $\gamma_1+\gamma_2=1$ subordinate to the open cover $\{A,B\}$.
Then for $(g,t)\in \supp (\gamma_2)$ we have

$$\distanza_{X}\left(\tilde{\mu}_{g^{-1}}^{G}\circ\tilde{\mu}_{t^{-1}}^{T}(x),x\right)\geq C\varepsilon$$
\noindent
for some constant $C>0$. Therefore, the map 

$$(g,t)\in P\mapsto \gamma_{2}(g,t)\chi_{\nu_{G}}(g^{-1})\Pi\left(\tilde{\mu}_{g^{-1}}^{G}\circ\tilde{\mu}_{t^{-1}}^{T}(x),x\right)$$ 
\noindent
is $\mathcal{C}^\infty$ because the singular support of $\Pi$ is included in the diagonal $X\times X$. As consequence the function:

$$ t\in T\mapsto \int_{G}\gamma_{2}(g,t)\chi_{\nu_{G}}(g^{-1})\Pi\left(\tilde{\mu}_{g^{-1}}^{G}\circ\tilde{\mu}_{t^{-1}}^{T}(x),x\right)dg $$
\noindent
is $\mathcal{C}^{\infty}$ and so its Fourier transform evaluated at $k\nu_{T}$ is rapidly decreasing for $k\rightarrow+\infty$, since by assumption $\nu_{T}\neq 0$.
We set $\gamma_{1}(g,\vartheta)=\sum_{j=1}^{N_{x}}\rho_{j}(g,\vartheta)$, with each $\rho_{j}$ supported in a neighborhood of $(g_{j},\vartheta_{j})$, and consider 

$$ \widetilde{\Pi}_{\nu_{G},k\nu_{T}}(x,x)\sim \sum_{j=1}^{N_{x}}\widetilde{\Pi}_{\nu_{G},k\nu_{T}}(x,x)^{(j)} $$
\noindent
where 

\begin{equation}
\label{arnold2}
\begin{multlined}[t][12.5cm]
\widetilde{\Pi}_{\nu_{G},k\nu_{T}}(x,x)^{(j)}\\
\sim \frac{d_{\nu_{G}}}{(2\pi)^{d_{T}}}\int_{G}\int_{(-\pi,\pi)^{d_{T}}}\rho_{j}(g,\vartheta)\chi_{\nu_{G}}(g^{-1})e^{-ik\nu_{T}\cdot\vartheta}\Pi\left(\tilde{\mu}_{g^{-1}}^{G}\circ\tilde{\mu}_{-\vartheta}^{T}(x),x\right)d\vartheta dg.
\end{multlined}                                             
\end{equation}

Let us now examine the asymptotics of each integrand separately.
On the support of $\gamma_1$, $\tilde{\mu}_{g^{-1}}^{G}\circ\tilde{\mu}_{t^{-1}}^{T}(x)$ is close to $x$, and therefore we may replace $\Pi$ by its representation as a Fourier integral, perhaps after disregarding a smoothing term which contributes negligibly to the asymptotics. After rescaling in $t$ we have:

\begin{equation}
\label{arnold3}
\begin{multlined}[t][12.5cm]
\widetilde{\Pi}_{\nu_{G},k\nu_{T}}(x,x)^{(j)}\\
\sim \frac{d_{\nu_{G}}}{(2\pi)^{d_{T}}}k\int_{G}\int_{(-\pi,\pi)^{d_{T}}}\int_{0}^{+\infty}\chi_{\nu_{G}}(g^{-1})e^{ik\Psi(x,t,g,\vartheta)}\\
\cdot s\left(\tilde{\mu}_{g^{-1}}^{G}\circ\tilde{\mu}_{-\vartheta}^{T}(x),x,kt\right)\gamma_{j}(g,\vartheta)d t d\vartheta dg
\end{multlined}
\end{equation}
\noindent
with 

\begin{equation}
\label{arnold4}
\begin{multlined}[t][12.5cm]
\Psi(x,t,g,\vartheta)=t\psi\left(\tilde{\mu}_{g^{-1}}^{G}\circ\tilde{\mu}_{-\vartheta}^{T}(x),x\right)-\nu_{T}\cdot\vartheta.
\end{multlined}
\end{equation}

Let us regard $\ref{arnold3}$ as an oscillatory integral with a complex phase $\Psi$ of positive type. Let us look for critical points of $\Psi$. 
We have that $\partial_{t}\Psi(x,t,g,\vartheta)=\psi\left(\tilde{\mu}_{g^{-1}}^{G}\circ\tilde{\mu}_{-\vartheta}^{T}(x),x\right)=0$ if and only if $\vartheta=\vartheta_{j}$ and $g=g_{j}$. We set $ \vartheta=\eta+\vartheta_{j}$ with $\eta\sim \mathbf{0}$.  In the neighborhood of $g_{j}$, we can write $g=g_{j}\exp_{G}{\xi}$, where $\xi\in \mathfrak{g}$ is close to the origin. With abuse of notation we shall write $\tilde{\mu}_{-\varsigma}^{G}\circ\tilde{\mu}_{-\vartheta}^{T}(x)=\tilde{\mu}_{-\xi}^{G}\circ\tilde{\mu}_{-\eta}^{T}(x)$.  Upon choosing orthonormal basis for the Lie algebras, we shall identify them with $\mathbb{R}^{d_{T}}$ and $\mathbb{R}^{d_{P}}$, respectively. We consider now $\left.\partial_{(\varsigma,\vartheta)}\Psi(x,t,\varsigma,\vartheta)\right|_{\xi=0,\eta=0}=t\Phi_{T}-\nu_{T}$. Thus in the new coordinates at any critical point $\xi=0$, $\eta=0$. The critical points are of the form:

$$P_{0}=(t_{0},(\varsigma_{0},\vartheta_{0}))=\left(\frac{\|\nu_{T}\|}{\|\Phi_{T}(m)\|},(\varsigma_{j},\vartheta_{j})\right).$$

Considering the second derivatives, $\partial_{tt}^{2}\Psi=0$, and using Heisenberg coordinates we write:

\begin{equation}
\label{noscalll}
\begin{multlined}[t][12.5cm]
\Psi(x,t,\xi,\eta)=t\psi\left(x+\left(\xi\Phi_{G}(m)+\eta\Phi_{T}(m)+ O(\|(\xi,\eta)\|^3),\right.\right.\\
\left.\left.,-\xi_{M}(m)-\eta_{M}(m)+ O(\|(\xi,\eta)\|^2)\right),x\right)-\nu_{T}\cdot\eta-\nu_{T}\cdot\vartheta_{j}.
\end{multlined}
\end{equation}

Here with abuse of language we have identified $\eta_{M}(m)$ and $\xi_{M}(m)$ with their representation in local coordinates.
We note that $m\in M_{0,\nu_{T}}$ so 

\begin{equation}
\label{noscalll2}
\begin{multlined}[t][12.5cm]
\Psi(x,t,\xi,\eta)=t\psi\left(x+\left(\eta\Phi_{T}(m)+ O(\|(\xi,\eta)\|^3),\right.\right.\\
\left.\left.,-\xi_{M}(m)-\eta_{M}(m)+ O(\|(\xi,\eta)\|^2)\right),x\right)-\nu_{T}\cdot\eta-\nu_{T}\cdot\vartheta_{j}\\
=-\nu_{T}\cdot\eta-\nu_{T}\cdot\vartheta_{j}+it\left\{\left[1-e^{i\eta\Phi_{T}(m)}\right]\right.\\ \left.+\frac{1}{2}\|\eta_{M}(m)+\xi_{M}(m)\|^{2}e^{i\eta\Phi_T(m)}
+ O(\|(\xi,\eta)\|^3)\right\}
\end{multlined}
\end{equation}
\noindent
and setting $\gamma=(\gamma_1,\cdots,\gamma_{d_P})=(\xi_{1},\cdots,\xi_{d_{G}},\eta_{1},\cdots,\eta_{d_{T}})$ and $\sigma=(\varsigma,\vartheta)$, we have $\sigma_{0}=(\varsigma_{0},\vartheta_{0})$ and

\begin{equation}
\label{noscalll2}
\begin{multlined}[t][12.5cm]
\left.\partial_{t\sigma}^{2}\Psi(x,t,\sigma)\right|_{(t_{0},\sigma_{0})}=\Phi_{T}(m).
\end{multlined}
\end{equation}

Thus, we have:

\begin{equation}
\label{noscalll3}
\begin{multlined}[t][12.5cm]
\left.\partial_{\sigma_{l}\sigma_{k}}^{2}\Psi(x,t,\sigma)\right|_{(t_{0},\sigma_{0})}=i\lambda_{\nu_{T}}\left[\Phi_{l}(m)\Phi_{k}(m)+\langle\gamma_{l},\gamma_{k}\rangle_{m}\right],
\end{multlined}
\end{equation}
\noindent
with $\lambda_{\nu_{T}}=t_{0}$, we are in the same case of $\cite{quarto}$ in the proof of Theorem 2, so we have that:

\begin{equation}
\label{superhessian}
\det{H(t_{0},\sigma_{0})}=i^{d_{P}+1}\lambda_{\nu_{T}}^{d_{P}-1}\|\Phi_{T}(m)\|^2\det{C(m)}.
\end{equation}

Where $C(m)$ is a scalar product on $\kernel{(\Phi_{P}(m))}$ and

\begin{equation}
\label{superhessian}
\det{\left[\frac{k}{2\pi i}H(t_{0},\sigma_{0})\right]}=\left(\frac{k}{2\pi}\right)^{d_{P}+1}\lambda_{\nu_{T}}^{d_{P}-1}\|\Phi_{T}(m)\|^2\mathcal{D}(m)^2.
\end{equation}

The principal term is:

\begin{equation}
\label{principal4}
\frac{d_{\nu_{G}}e^{-ik\vartheta_{j}\cdot\nu_{T}}2^{d_{G}/2}\chi_{\nu_{G}}(g_{j}^{-1})}{\left(\sqrt{2}\pi\right)^{d_{T}-1}}\cdot\left(\frac{\|\nu_{T}\|k}{\pi}\right)^{d_{M}+\frac{1-d_{P}}{2}}\cdot\frac{1}{\mathcal{D}(m)\|\Phi_{T}(m)\|^{d_{M}+1+\frac{1-d_{P}}{2}}}.
\end{equation}

This complete the proof of the Theorem.
\hfill $\Box$

\section{Proof of Corollary $\bf{~\ref{cor:quarto}}$}

$\Proof.$

We start considering the dimension of $H(X)_{\nu_{G},k\nu_{T}}$:

$$\dimensione{\left(H(X)_{\nu_{G},k\nu_{T}}\right)}=\int_{X}\widetilde{\Pi}_{\nu_{G},k\nu_{T}}(x,x)dV_{X}(x).$$

Now let us observe that $\widetilde{\Pi}_{\nu_{G},k\nu_{T}}(x,x)$ is naturally $S^{1}$-invariant as a function of $x$, and therefore descends to a function
on $M$, that we shall denote by $\widetilde{\Pi}_{\nu_{G},k\nu_{T}}(m,m)$ with abuse of language. Thus by integrating first along the fibers the previous integral may be naturally interpreted as an integral over $M$, that we shall write in the form:

$$\dimensione{\left(H(X)_{\nu_{G},k\nu_{T}}\right)}=\int_{M}\widetilde{\Pi}_{\nu_{G},k\nu_{T}}(m,m)dV_{M}(m).$$

Now by the above $\widetilde{\Pi}_{\nu_{G},k\nu_{T}}(m,m)$ is rapidly decreasing away from a shrinking neighborhood of $M_{0,\nu_{T}}$. So, using a smoothly varying system of adapted coordinates centered at points $m\in M_{0,\nu_{T}}$, we can locally parametrize a neighborhood $U$ of $M_{0,\nu_{T}}$ in the form $m+\vv$, where $m\in M_{0,\nu_{T}}$ and 
$\vv\in N_{m}$. This parametrization is only valid locally in $m$, since we may not expect to find a single $\mathcal{C}^\infty$ family of adapted coordinates $\gamma_{m}$ ($m\in M_{0,\nu_{T}}$). Hence to make this argument strictly rigorous we should introduce a partition of unity on $M_{0,\nu_{T}}$ subordinate to an appropriate open cover. However, we shall simplify notation and leave this point implicit.

\begin{equation*}
\dimensione{\left(H(X)_{\nu_{G},k\nu_{T}}\right)}=\int_{M_{0,\nu_{T}}}\int_{\mathbb{R}^{d_{P}-1}}\widetilde{\Pi}_{\nu_{G},k\nu_{T}}\left(m+\vv,m+\vv\right)d\vv dV_{M}(m).
\end{equation*}

In view of Theorem $~\ref{teo:secondo}$ the asymptotics of the previous integral are unchanged, if the integrand is multiplied by a cut-off of the form $\varrho \big(k^{\frac{7}{18}}\|\vv\|\big)$, where $\varrho\in \mathcal{C}^\infty_{0}(\mathbb{R})$ is identically equal to $1$ in some neighborhood of $0$.

\begin{equation*}
\dimensione{\left(H(X)_{\nu_{G},k\nu_{T}}\right)}=\int_{M_{0,\nu_{T}}}\int_{\mathbb{R}^{d_{P}-1}}\widetilde{\Pi}_{\nu_{G},k\nu_{T}}\left(m+\vv,m+\vv\right)\varrho \big(k^{\frac{7}{18}}\|\vv\|\big)d\vv dV_{M}(m).
\end{equation*}

Let us now operate the rescaling $\vv=\frac{\uu}{\sqrt{k}}$. We can now make use of the asymptotic expansion in Theorem $~\ref{teo:diagonal}$, with $\uu=\uu_{t}$ (that is, $\uu_{v}=\uu_{h}=0$). We obtain:

\begin{equation}
\label{telufficio5731}
\begin{multlined}[t][12.5cm]
\dimensione{\left(H(X)_{\nu_{G},k\nu_{T}}\right)}=k^{-\frac{d_{P}-1}{2}}\int_{M_{0,\nu_{T}}}\int_{\mathbb{R}^{d_{P}-1}}\widetilde{\Pi}_{\nu_{G},k\nu_{T}}\left(m+\frac{\uu}{\sqrt{k}},m+\frac{\uu}{\sqrt{k}}\right)\varrho \big(k^{-\frac{1}{9}}\|\uu\|\big)d\uu dV_{M}\\
=k^{-\frac{d_{P}-1}{2}}\cdot\int_{M_{0,\nu_{T}}}\frac{2^{\frac{d_{G}}{2}}d_{\nu_{G}}^{2}}{(\sqrt{2})^{d_{T}-1}\pi^{d_{T}-1}}\left(\frac{\|\nu_{T}\|k}{\pi}\right)^{d_{M}-\frac{d_{P}-1}{2}}\cdot\\
\cdot\int_{\mathbb{R}^{d_{P}-1}}\frac{1}{\|\Phi_{T}\|^{d_{M}-\frac{d_{P}-1}{2}+1}\det{C(m)}}e^{-\lambda_{\nu_{T}}2\|\uu\|^2}\varrho \big(k^{-\frac{1}{9}}\|\uu\|\big)d\uu dV_{M}(m) + \cdots,
\end{multlined}
\end{equation}
\noindent
where $d_{P}=d_{G}+d_{T}$ and the dots denote lower order terms.
Now we evaluate the Gaussian integral, let us operate the change of variables $q=\sqrt{2\lambda_{\nu_{T}}}\uu$, we have that:

\begin{equation}
\label{boiaboiagausss}
\begin{multlined}[t][12.5cm]
\int_{\mathbb{R}^{d_{P}-1}} e^{-2\lambda_{\nu_{T}}\|\uu\|^{2}}d\uu=\frac{1}{\lambda_{\nu_{T}}^{\frac{d_{P}-1}{2}}(\sqrt{2})^{d_{P}-1}}\int_{\mathbb{R}^{d_{P}-1}} e^{-\|q\|^{2}}dq\\
=\frac{\|\Phi_{T}(m)\|^{\frac{d_{P}-1}{2}}\pi^{\frac{d_{P}-1}{2}}}{(\sqrt{2})^{d_{P}-1}\|\nu_{T}\|^{\frac{d_{P}-1}{2}}}
\end{multlined}
\end{equation}
\noindent
and substituting in $(\ref{telufficio5731})$ we obtain the following expression:

\begin{equation*}
\begin{multlined}[t][12.5cm]
\dimensione{\left(H(X)_{\nu_{G},k\nu_{T}}\right)}=\frac{d_{\nu_{G}}^{2}}{2^{d_{T}-1}\pi^{d_{T}-1}}\left(\frac{\|\nu_{T}\|k}{\pi}\right)^{d_{M}-d_{P}+1}\cdot\\
\cdot\int_{M_{0,\nu_{T}}}\frac{\|\Phi_{T}(m)\|^{-d_{M}+d_{P}-2}}{\det{C(m)}}dV_{M}(m)+ \cdots.
\end{multlined}
\end{equation*}

The proof is complete.
\hfill $\Box$

\section{Proof of Proposition $\bf{~\ref{prop:complement}}$}

$\Proof.$

We assume that $\distanza_{X}{(y,p\cdot x)}\geq Dk^{\varepsilon-1/2}$, for every $p\in P$. The method consist to use iteratively integration by parts to deduce the rapidly decreasing behavior of the kernel. We start following the previous situations using the standard representation as Fourier integral operator. So, by performing the change of variables $t\mapsto kt$, we obtain the following expression for the Szeg\"{o} kernel:  

\begin{equation}
\label{geometriainbicocca2015}
\begin{multlined}[t][12.5cm]
\widetilde{\Pi}_{\nu_{G},k\nu_{T}}\left(x,y\right)=\frac{kd_{\nu_{G}}}{(2\pi)^{d_{T}}}\int_{G}\int_{T}\int_{0}^{+\infty}\chi_{\nu_{G}}(g^{-1})e^{ik\Psi(t,\vartheta,x)}\cdot\\
\cdot s\left(\tilde{\mu}^{G}_{g^{-1}}\circ\tilde{\mu}^{T}_{-\vartheta}\left(x\right),y,kt\right)dtdgd\vartheta,
\end{multlined}
\end{equation}
\noindent
where 

\begin{equation}
\label{geometriainbicocca20152}
\Psi(t,\vartheta,x)=t\psi\left(\tilde{\mu}^{G}_{g^{-1}}\circ\tilde{\mu}^{T}_{-\vartheta}\left(x\right),y\right)-\nu_{T}\cdot \vartheta
\end{equation}
\noindent
is the phase of the oscillatory integral. First we observe that:

$$ \|\partial_{\vartheta}\Psi\|=\|t\psi\left(\tilde{\mu}^{G}_{g^{-1}}\circ\tilde{\mu}^{T}_{-\vartheta}\left(x\right),y\right)-\nu_{T}\|\geq C$$
\noindent
for $0<t\ll 1$ and $C>0$. In a similar way 

$$ \|\partial_{\vartheta}\Psi\|\geq C_{1}t $$
\noindent 
for $t\gg 0$ and $C_{1}>0$. Using integration by parts in $d\vartheta$, the asymptotics for $k\rightarrow +\infty$ is unchanged. We multiply the integrand by $\gamma(t)$, where $\gamma\in \mathcal{C}^{\infty}_{0}\left(\frac{1}{2D},2D\right)$, $\gamma\equiv 1$ on $\left(\frac{1}{D},D\right)$ and $\gamma \equiv 0$ outside of $\left(\frac{1}{2D},2D\right)$. The integral now is compactly supported in $dt$. Taking the partial derivative respect $t$, we deduce that:

$$\partial_{t}\Psi(t,\vartheta,x)=\psi\left(\tilde{\mu}^{G}_{g^{-1}}\circ\tilde{\mu}^{T}_{-\vartheta}\left(x\right),y\right)$$
\noindent
and by the assumption we find that:

$$\left\|\partial_{t}\Psi(t,\vartheta,x)\right\|=\left|\psi\left(\tilde{\mu}^{G}_{g^{-1}}\circ\tilde{\mu}^{T}_{-\vartheta}\left(x\right),y\right)\right|\geq \Im{\psi\left(\tilde{\mu}^{G}_{g^{-1}}\circ\tilde{\mu}^{T}_{-\vartheta}\left(x\right),y\right)}\geq D'k^{2\varepsilon-1}.$$

Now we introduce the differential operator:

$$ L_{t}=\left[\psi\left(\tilde{\mu}^{G}_{g^{-1}}\circ\tilde{\mu}^{T}_{-\vartheta}\left(x\right),y\right)\right]^{-1}\partial_{t} $$
\noindent
and observing that $e^{ik\Psi}=-\frac{i}{k}L_{t}\left(e^{ik\Psi}\right)$ we can apply iteratively the integration by parts. So step by step we obtain a factor of order $O(k^{-2N\varepsilon})$.

The proof is complete.
\hfill $\Box$

\section{Proof of Theorem $\bf{~\ref{teo:app2}}$}

$\Proof.$

In view of the equality on the first line of $(\ref{toeplitz})$ we observe that $1)$ follows immediately from the point $1)$ of the main Theorem.
Let us now consider the proof of $2)$.

Let $f\in C^\infty(M)$, we consider the associated Toeplitz operator:

\begin{equation}
\label{Toeplitz1}
\begin{multlined}[t][12.5cm]
T_{\nu_{G},k\nu_{T}}[f]\left(x+\frac{\vv}{\sqrt{k}},x+\frac{\vv}{\sqrt{k}}\right)\\=\int_{X}\widetilde{\Pi}_{\nu_{G},k\nu_{T}}\left(x+\frac{\vv}{\sqrt{k}},y\right)f(y)\widetilde{\Pi}_{\nu_{G},k\nu_{T}}\left(y,x+\frac{\vv}{\sqrt{k}}\right)dV_{X}(y)
\end{multlined}
\end{equation}
\noindent
with $x\in X_{0,\nu_{T}}$, where $f(y)=f(\pi(y))$. Now in view of Proposition $\ref{prop:complement}$ only a shrinking neighborhood of the orbit $P\cdot x$ contributes non-negligibly to the asymptotics. Therefore, the asymptotics are unchanged if the integrand in $(\ref{Toeplitz1})$ is multiplied by a cut-off function $\varrho_{k}(y)$, where $\varrho_{k}=1$ for $\distanza_{X}(y,P\cdot x)\leq Dk^{\delta-1/2}$ (for example concretely $\delta$ equal $1/9$) and $\varrho_{k}=0$ for $\distanza_{X}{(y,P\cdot x)}\geq 2Dk^{\delta-1/2}$. We shall make a more explicit choice of $\varrho_{k}$ below.

Let $x+(\theta,\vv)$ be a system of Heisenberg local coordinates on $X$ centered at $x$. This determines for every $p\in P$ a system of HLC centered at $p\cdot x$, by setting 

$$p\cdot x+(\theta,\vv)= p\cdot \big(x+(\theta,\vv)\big).$$

In this manner we have a unique smoothly varying family of HLC systems centered at points of $P\cdot x$, and identifications $T_{p\cdot x}X\cong T_{x}X\cong 
\mathbb{R}\times\mathbb{R}^{2d_{M}}$, $T_{p\cdot m}M\cong T_{m}M\cong \mathbb{R}\times\mathbb{R}^{2d_{M}}$. Furthermore, the action of $P$ preserves the contact and CR structures of $M$, and the decomposition of the tangent spaces in $h$-, $v$-, and $t$-components. This means that the corresponding decomposition is preserved under the identification $T_{p\cdot m}M\cong T_{m}M$. If $y\in P\cdot x$, let $N^{P}_{y}$ be the normal space to $P\cdot x$ in $X$ at $y$; then we have natural unitary isomorphisms $N^{P}_{y}\cong N^{P}_{x}$. 
With this identification implicit and some abuse of language, we may then parametrize a suitably small open neighborhood of $P\cdot x$ by the map 

$$(p,n)\in P\times N^{P}_{x}\mapsto p\cdot x+n.$$

We set $y=p\cdot \left(x+n\right)$, where $p=(g,t)\in P$ and $n$ is a tangent vector normal to the orbit. 
For simplicity we suppose that the stabilizator of $x$ in $P$ is trivial. We have a diffeomorphism $P\times N(\varepsilon)\rightarrow X'$, with $X'$ a $\varepsilon$-tubular neighborhood of $P\cdot x$ in $X$ and $N(\varepsilon)$ is a ball of radious $\varepsilon$ in the normal space of the orbit in $x$ (a real vector space of dimension $2d_{M}+1-d_{P}$). This diffeomorphism doesn't preserve the volume form. Thus in coordinates $(p,n)\in P\times N_{x}^{P}$ (where $N_{x}^{P}$ denote the normal space to $P\cdot x$ in $x$):

\begin{equation}
\label{newcoordinates}
dV_{X}(y)=\mathfrak{D}(p,n)dV_{P}(p)dn,
\end{equation}
\noindent
where $dV_{P}(p)=\frac{d\vartheta}{(2\pi)^{d_{T}}}dV_{G}(g)$ is the Haar measure of $P$, $dn=d\mathcal{L}(n)$ the Lebesgue measure on  $N_{x}^{P}$ (unitarily identified with $\mathbb{R}^{2d_{M}+1-d_{P}}$) and $\mathfrak{D}(p,0)=\mathfrak{R}_{x}(p)$ with $\mathfrak{R}_{x}:P\rightarrow \mathbb{R}_{>0}$ a distortive function defined as follow. 
Let $\mathcal{B}_{0}$ an orthonormal basis of $\mathfrak{p}$ and let $\valutation_{x}:\mathfrak{p}\rightarrow T_{x}X$ the valutation map. Let $D_{x}(p)$ the matrix associated to $\valutation_{x}^{\vee}(g_{X}):\mathfrak{p}\times \mathfrak{p}\rightarrow \mathbb{R}$ respect to $\mathcal{B}_{0}$ and 

\begin{equation}
\label{distortionf}
\mathfrak{R}_{x}(p)=\sqrt{\det{D_{x}(p)}},
\end{equation}
\noindent
here $g_{X}$ is the Riemannian metric on $X$ and $\valutation_{x}^{\vee}(g_{X})$ is the pull back of such metric to $\mathfrak{p}=\lie{(P)}$ using the valutation 
$\valutation_{x}:\gamma \mapsto \gamma_{X}(x)$.
We observe that $\mathfrak{R}_{x}(p)$ is constant along the orbit, because $P$ acts by isometries.
We put $\mathfrak{R}_{x}(p)=r_{x}$. 
We obtain that:

\begin{equation}
\label{Toeplitz2}
\begin{multlined}[t][12.5cm]
T_{\nu_{G},k\nu_{T}}[f]\left(x+\frac{\vv}{\sqrt{k}},x+\frac{\vv}{\sqrt{k}}\right)\\
\\
=\int_{X}\widetilde{\Pi}_{\nu_{G},k\nu_{T}}\left(x+\frac{\vv}{\sqrt{k}},y\right)f(y)\widetilde{\Pi}_{\nu_{G},k\nu_{T}}\left(y,x+\frac{\vv}{\sqrt{k}}\right)\varrho_{k}(y)dV_{X}(y)\\
=\int_{P}\int_{N_{x}}\widetilde{\Pi}_{\nu_{G},k\nu_{T}}\left(x+\frac{\vv}{\sqrt{k}},p\left(x+n\right)\right)f\left(p\cdot\left(x+n\right)\right)\cdot\\
\cdot\widetilde{\Pi}_{\nu_{G},k\nu_{T}}\left(p\left(x+n\right),x+\frac{\vv}{\sqrt{k}}\right)\mathfrak{D}(p,n)\varrho_{k}\left(p\cdot\left(x+n\right)\right)dV_{P}(p)dn.
\end{multlined}
\end{equation}

Rescaling $n\mapsto \frac{n}{\sqrt{k}}$, we observe that $\rk{(N)}=\dimensione{X}-\dimensione{P}=2d_{M}+1-d_{P}=2\left[d_{M}+\frac{1-d_{P}}{2}\right]$ from which we have that 
$dn \rightarrow k^{-\left[d_{M}+\frac{1-d_{P}}{2}\right]}dn$ and we obtain that:

\begin{equation}
\label{Toeplitz3}
\begin{multlined}[t][12.5cm]
T_{\nu_{G},k\nu_{T}}[f]\left(x+\frac{\vv}{\sqrt{k}},x+\frac{\vv}{\sqrt{k}}\right)\\ \\=k^{-\left[d_{M}+\frac{1-d_{P}}{2}\right]}\int_{X}\widetilde{\Pi}_{\nu_{G},k\nu_{T}}\left(x+\frac{\vv}{\sqrt{k}},y\right)f(y)\widetilde{\Pi}_{\nu_{G},k\nu_{T}}\left(y,x+\frac{\vv}{\sqrt{k}}\right)\varrho_{k}(y)dV_{X}(y)\\
=k^{-\left[d_{M}+\frac{1-d_{P}}{2}\right]}\int_{P}\int_{N_{x}}\widetilde{\Pi}_{\nu_{G},k\nu_{T}}\left(x+\frac{\vv_{1}}{\sqrt{k}},p\left(x+\frac{n}{\sqrt{k}}\right)\right)f\left(p\left(x+\frac{n}{\sqrt{k}}\right)\right)\cdot\\
\cdot\widetilde{\Pi}_{\nu_{G},k\nu_{T}}\left(p\left(x+\frac{n}{\sqrt{k}}\right),x+\frac{\vv}{\sqrt{k}}\right)\mathfrak{D}\left(p,\frac{n}{\sqrt{k}}\right)\varrho_{k}\left(p\cdot \left(x+\frac{n}{\sqrt{k}}\right)\right)dV_{P}(p)dn.
\end{multlined}
\end{equation}

Here $\varrho_k(p\cdot x+n)=\varrho_k\left(p\cdot x+k^{1/2-\varepsilon}n\right)$ and $\varepsilon=1/9$. $n\in N_{x}^{P}$ by construction and $v\in N^{P}_{x}$ by assumption. We have that: 

$$N_{x}^{P}=[\valutation_{x}(\kernel{(\Phi_{P}(m))})\oplus\spaned{(\eta_{X}(x))}]^{\perp}.$$

Let $\eta\in \mathfrak{t}$ be the unique element such that
$$
\eta\in \ker\Phi_{T}(m)^\perp,\,\,\,\,\langle \Phi_{T}(m),\eta\rangle=\|\Phi_{T}(m)\|.
$$
In particular we have $\eta\in \mathfrak{t}$, and $\eta$ has unit norm. Then
 
\begin{equation}
\label{HLCinx}
\eta_{X}(x)=(\|\Phi_{P}(m)\|,-\eta_{M}(m)).
\end{equation}

So in terms of the isomorphism $T_{x}X\cong \mathbb{R}\times T_{m}M$ we have:

\begin{equation*}
p_{X}(x)=(\{0\}\times\valutation_{m}(\kernel{(\Phi_{P}(m))}))\oplus\spaned_{\mathbb{R}}{((\|\Phi(m)\|,-\eta_{M}(m)))},
\end{equation*}
\noindent
and so in the notation of $(\ref{labb})$ and $(\ref{antoniosalierii})$

\begin{equation*}
\begin{multlined}[t][12.5cm]
p_{X}(x)^{\perp}=(\{0\}\times V_{m})^{\perp}\cap\spaned_{\mathbb{R}}{(\|\Phi(m)\|,-\eta_{M}(m))}^{\perp}\\
=\left[\mathbb{R}\times(H_{m}\oplus N_{m})\right]\cap[\spaned_{\mathbb{R}}{((\|\Phi(m)\|,-\eta_{M}(m)))}]^{\perp},
\end{multlined}
\end{equation*}
\noindent
with $V_{m}=\valutation_{m}(\kernel{(\Phi_{P}(m))})$, $N_{m}=J_{m}(V_{m})$ and $H_{m}=(V_{m}\oplus N_{m})^{\perp}$. 
We have that $\eta_{Mt}(m)=0$, because $M_{0,\nu_{T}}$ is $P$-invariant. Thus if $(\lambda,\hh+\ttt)\in \mathbb{R}\times(H_{m}\oplus N_{m})$ we have:

\begin{equation*}
g_{X}((\lambda,\hh+\ttt),(\|\Phi(m)\|,-\eta_{M}(m)))=\lambda\|\Phi_{P}(m)\|-g_{X}(\hh,\eta_{Mh}(m)).
\end{equation*}

Denoting $\Phi=\Phi_{P}(m)$ and $\eta_{M}=\eta_{M}(m)$ we have:

\begin{equation}
\label{normalspaaace}
\begin{multlined}[t][12.5cm]
N_{x}^{P}=\left\{(\lambda,\hh+\ttt)\in\mathbb{R}\times(H_{m}\oplus N_{m}): \lambda=\frac{g_{X}(\hh,\eta_{Mh})}{\|\Phi\|}\right\}\\
=[\{0\}\times N_{m}]\oplus\left\{\left(\frac{1}{\|\Phi\|}g_{X}(\hh,\eta_{Mh}),\hh\right):\hh\in H_{m}\right\}.
\end{multlined}
\end{equation}

We write $v=n_{1}=\left(\frac{1}{\|\Phi\|}g_{X}(\hh_{1},\eta_{Mh}),\hh_{1}+\ttt_{1}\right)$, $n=\left(\frac{1}{\|\Phi\|}g_{X}(\hh,\eta_{Mh}),\hh+\ttt\right)$. Recalling that $\Pi(x_1,x_2)=\overline{\Pi(x_2,x_1)}$ we obtain:

\begin{equation}
\label{Toeplitz4}
\begin{multlined}[t][12.5cm]
T_{\nu_{G},k\nu_{T}}[f]\left(x+\frac{n_{1}}{\sqrt{k}},x+\frac{n_{1}}{\sqrt{k}}\right)=\\
k^{-\left[d_{M}+\frac{1-d_{P}}{2}\right]}\int_{P}\int_{N_{x}}\left|\widetilde{\Pi}_{\nu_{G},k\nu_{T}}\left(x+\frac{n_{1}}{\sqrt{k}},p\left(x+\frac{n}{\sqrt{k}}\right)\right)\right|^{2}f\left(p\left(x+\frac{n}{\sqrt{k}}\right)\right)\cdot\\
\cdot\mathfrak{D}\left(p,\frac{n}{\sqrt{k}}\right)\varrho_{k}\left(p\cdot \left(x+\frac{n}{\sqrt{k}}\right)\right)dV_{P}(p)dn.
\end{multlined}
\end{equation}

Let us now make use of  the asymptotic expansion of $\widetilde{\Pi}_{\nu_{G},k\nu_{T}}$ from point 3) of the Theorem $\ref{teo:secondo}$ and, using the Taylor expansion for $f\left(p\cdot x+\frac{n}{\sqrt{k}}\right)$ and for $\mathfrak{D}$

\begin{equation}
\label{Toeplitz6}
\begin{multlined}[t][12.5cm]
T_{\nu_{G},k\nu_{T}}[f]\left(x+\frac{n_{1}}{\sqrt{k}},x+\frac{n_{1}}{\sqrt{k}}\right)\sim k^{\left[d_{M}+\frac{1-d_{P}}{2}\right]}r_{x}\cdot f\left(p\cdot x\right)C(m,\nu_{P})\\ 
\int_{P}\int_{N_{x}}\left|\chi_{\nu_{P}}\left(p\right)\right|^2\cdot e^{H(n_{1},n)+H(n,n_{1})}\varrho_{k}\left(p\cdot \left(x+\frac{n}{\sqrt{k}}\right)\right)dV_{P}(p)dn(1+\cdots)
\end{multlined}
\end{equation}
\noindent
where we have set:

$$ C(m,\nu_{P})=\left[\frac{d_{\nu_{G}}2^{\frac{d_{G}}{2}}}{(\sqrt{2}\pi)^{d_{T}-1}} 
\cdot\frac{\left(\frac{\|\nu_{T}\|}{\pi}\right)^{d_{M}+\frac{1-d_{P}}{2}}}{\mathcal{D}(m)\|\Phi_{T}\|^{d_{M}+1+\frac{1-d_{P}}{2}}}\right]^2 $$
\noindent
and the dots stand for terms of less degree. 

Here the exponent is as follows: let 
$$n=\vv+\ttt+\hh,\,\,\,\,\, n_1=\vv_{1}+\ttt_{1}+\hh_{1}$$
be the decomposition as in $(\ref{labb})$ and $(\ref{antoniosalierii})$. We call $R_{x}=N_{x}^{P}$ and $r,r_{1}$ in place of $n_{1},n$. Then

\begin{equation}
\label{exxponent1}
\begin{multlined}[t][12.5cm]
H(r_{1},r)+H(r,r_{1})=-2\lambda_{\nu_{T}}\left(\|\ttt_{1}\|^2+\|\ttt\|^2\right)-\lambda_{\nu_{T}}\left\|\hh_{1}-\hh-\frac{(\hh-\hh_{1},\eta_{Mh})}{\|\Phi_{T}(m)\|}\eta_{Mh}\right\|^2.
\end{multlined}
\end{equation}

We evaluate the Gaussian integral $\int_{R_{x}}e^{H(r_{1},r)+H(r,r_{1})}dr$. So we have that:

\begin{equation}
\label{Gaussian11}
\int_{R_{x}}e^{H(r_{1},r)+H(r,r_{1})}dr=e^{-2\lambda_{\nu_{T}}\|\ttt_{1}\|^2}\int_{N_{m}}e^{-2\lambda_{\nu_{T}}\|\ttt\|^{2}}d\ttt\cdot\int_{H_{m}}e^{-\lambda_{\nu_{T}}\left\|\hh_{1}-\hh-\frac{(\hh-\hh_{1},\eta_{Mh})}{\|\Phi_{T}(m)\|}\eta_{Mh}\right\|^2}d\hh.
\end{equation}

Let $d_{N_{m}}=d_{P}-1$ the dimension of $N_{m}$. Let us first consider the first Gaussian integral in $(\ref{Gaussian11})$.

\begin{equation}
\label{Gaussian22}
\int_{N_{m}}e^{-2\lambda_{\nu_{T}}\|\ttt\|^{2}}d\ttt=\frac{1}{\left(2\sqrt{\lambda_{\nu_{T}}}\right)^{d_{N_{m}}}}\int_{N_{m}}e^{-\frac{1}{2}\|s\|^{2}}ds=\left(\frac{\pi}{2\lambda_{\nu_{T}}}\right)^{(d_P-1)/2}.
\end{equation}

To compute the second Gaussian integral in $(\ref{Gaussian11})$, let us operate the change of variable 

\begin{equation}
\label{samsung44}
-\hh_{1}+\hh+\frac{(\hh-\hh_{1},\eta_{Mh})}{\|\Phi_{T}(m)\|}\eta_{Mh}=\ww.
\end{equation}

We differentiate the previous expression $d\hh=\det{\left(\frac{\partial \ww}{\partial \hh}\right)}^{-1}d\ww$ and in order to determine $\det{\left(\frac{\partial \ww}{\partial \hh}\right)}$ we observe that:

$$\det{\left(\frac{\partial \ww}{\partial \hh}\right)}=1+\frac{\|\eta_{Mh}\|^2}{\|\Phi_{T}\|}.$$

So we have:

\begin{equation}
\label{Gaussian33}
\begin{multlined}[t][12.5cm]
\int_{H_{m}}e^{-\lambda_{\nu_{T}}\left\|\hh_{1}-\hh-\frac{(\hh-\hh_{1},\eta_{Mh})}{\|\Phi_{T}(m)\|}\eta_{Mh}\right\|^2}d\hh\\
=\frac{\|\Phi_{T}\|}{\|\Phi_{T}\|+\|\eta_{Mh}\|^2}\int_{H_{m}}e^{-\frac{1}{2}\|s\|^2}\frac{ds}{\left(2\lambda_{\nu}\right)^{\frac{d_{H_{m}}}{2}}}=\frac{\|\Phi_{T}\|}{\|\Phi_{T}\|+\|\eta_{Mh}\|^2}\left(\frac{\pi}{\lambda}\right)^{d_{H_{m}}/2},
\end{multlined}
\end{equation}
\noindent
with $d_{H_{m}}=2(d_{M}+1-d_{P})$ the dimension of $H_{m}$.
Inserting $(\ref{Gaussian33})$ and $(\ref{Gaussian22})$ in $(\ref{Gaussian11})$,
we obtain

\begin{equation}
\label{Gaussianfinalll}
\begin{multlined}[t][12.5cm]
\int_{R_{x}}e^{H(r_{1},r)+H(r,r_{1})}dr\\
=e^{-2\lambda_{\nu_{T}}\|\ttt_{1}\|^2}\cdot\frac{\|\Phi_{T}\|}{\|\Phi_{T}\|+\|\eta_{Mh}\|^2}\left(\frac{\pi}{\lambda_{\nu_{T}}}\right)^{\frac{d_{H_{m}}+d_{N_{m}}}{2}}2^{-\frac{d_{N_{m}}}{2}}.
\end{multlined}
\end{equation}

Let us now insert this in $(\ref{Toeplitz6})$. We obtain a leading term depending on $r_x$ and $f(x)$. We can then determine $r_x$ using that for $f=1$ this must reduce to the leading term in 2) of Theorem $\ref{teo:secondo}$. Thus, recalling that $\int_{P}\left|\chi_{\nu_{P}}\left(p\right)\right|^2dV_{P}(p)=1$ (see $\cite{ventisettesimo}$ Theorem 4.11) and noting that $d_{H_{m}}+d_{N_{m}}=2d_{M}+1-d_P$ we obtain that the leading order term in $(\ref{Toeplitz6})$ is given by:

\begin{equation}
\label{costanteperduta}
\begin{multlined}[t][12.5cm]
k^{\left[d_{M}+\frac{1-d_{P}}{2}\right]}\cdot C(m,\nu_{P})\cdot\left(\frac{\pi}{\lambda_{\nu_{T}}}\right)^{\frac{d_{H_{m}}+d_{N_{m}}}{2}}\cdot 2^{-\frac{d_{N_{m}}}{2}}\\
\cdot e^{-2\lambda_{\nu_{T}}\|\ttt_{1}\|^2}\cdot\frac{\|\Phi_{T}\|}{\|\Phi_{T}\|+\|\eta_{Mh}\|^2} \cdot r_{x}=\frac{1}{(\sqrt{2}\pi)^{d_{T}-1}}d_{\nu_{G}}2^{\frac{d_{G}}{2}}\left(\frac{k}{\pi}\|\nu_{T}\|\right)^{d_{M}-\frac{d_{P}}{2}+\frac{1}{2}}e^{-2\lambda_{\nu_{T}}\|\ttt_{1}\|^2}\cdot\\
\cdot \frac{1}{\mathcal{D}(m)}\cdot\frac{1}{\|\Phi_{T}\|^{d_{M}+1-\frac{d_{P}}{2}+\frac{1}{2}}},
\end{multlined}
\end{equation}
\noindent
and we find that 

\begin{equation}
\label{costanteperduta2}
\begin{multlined}[t][12.5cm]
r_{x}=\frac{\pi^{d_{T}-1}\mathcal{D}(m)(\|\Phi_{T}\|+\|\eta_{Mh}\|^2)}{d_{\nu_{G}}\left(\sqrt{2}\right)^{-d_{G}-d_{T}+1-d_{N_{m}}}}.
\end{multlined}
\end{equation}

The leading term become:

\begin{equation}
\label{Toeplitz7}
\begin{multlined}[t][12.5cm]
\frac{1}{(\sqrt{2}\pi)^{d_{T}-1}}d_{\nu_{G}}2^{\frac{d_{G}}{2}}\left(\frac{k}{\pi}\|\nu_{T}\|\right)^{d_{M}-\frac{d_{P}}{2}+\frac{1}{2}}f(m)e^{-2\lambda_{\nu_{T}}\|\ttt_{1}\|^2}\cdot\\
\cdot \frac{1}{\mathcal{D}(m)}\cdot\frac{1}{\|\Phi_{T}\|^{d_{M}+1-\frac{d_{P}}{2}+\frac{1}{2}}}.
\end{multlined}
\end{equation}

This complete the proof of $2)$ and of the Theorem.
\hfill $\Box$

\section{Proof of Corollary $\bf{~\ref{cor:quinto}}$}

$\Proof.$

We start considering the trace of $T_{\nu_{G},k\nu_{T}}[f]$:

$$\mathfrak{T}\left(T_{\nu_{G},k\nu_{T}}[f]\right)=\int_{X}T_{\nu_{G},k\nu_{T}}[f](x,x)dV_{X}(x).$$

Now we observe that $T_{\nu_{G},k\nu_{T}}(x,x)$ is rapidly decreasing away from a shrinking neighborhood of $X_{0,\nu_{T}}$. So, using a smoothly varying system of adapted coordinates centered at points $x\in X_{0,\nu_{T}}$, we can locally parametrize a neighborhood $U$ of $X_{0,\nu_{T}}$ in the form $x+\ttt$, where $x\in X_{0,\nu_{T}}$ and 
$\ttt\in N_{m}$. This parametrization is only valid locally in $x$. We introduce a partition of unity on $X_{0,\nu_{T}}$ subordinate to an appropriate open cover and we simplify the notation leaving this point implicit.

\begin{equation*}
\mathfrak{T}\left(T_{\nu_{G},k\nu_{T}}[f]\right)=\int_{X_{0,\nu_{T}}}\int_{\mathbb{R}^{d_{P}-1}}T_{\nu_{G},k\nu_{T}}[f]\left(x+\ttt,x+\ttt\right)d\ttt dV_{X}(x).
\end{equation*}

In view of Theorem $~\ref{teo:app2}$ the asymptotics of the previous integral are unchanged, if the integrand is multiplied by a cut-off of the form $\varrho \big(k^{\frac{7}{18}}\|\ttt\|\big)$, where $\varrho\in \mathcal{C}^\infty_{0}(\mathbb{R})$ is identically equal to $1$ in some neighborhood of $0$.

\begin{equation*}
\mathfrak{T}\left(T_{\nu_{G},k\nu_{T}}[f]\right)=\int_{X_{0,\nu_{T}}}\int_{\mathbb{R}^{d_{P}-1}}T_{\nu_{G},k\nu_{T}}\left(m+\vv,m+\vv\right)\varrho \big(k^{\frac{7}{18}}\|\ttt\|\big)d\ttt dV_{X}(x).
\end{equation*}

Let us now operate the rescaling $\ttt=\frac{\uu}{\sqrt{k}}$. We can now make use of the asymptotic expansion in Theorem $~\ref{teo:app2}$, with $n_{1}=\uu$. We obtain:

\begin{equation}
\label{telufficio57314}
\begin{multlined}[t][12.5cm]
\mathfrak{T}\left(T_{\nu_{G},k\nu_{T}}[f]\right)=k^{-\frac{d_{P}-1}{2}}\int_{X_{0,\nu_{T}}}\int_{\mathbb{R}^{d_{P}-1}}T_{\nu_{G},k\nu_{T}}\left(x+\frac{\uu}{\sqrt{k}},x+\frac{\uu}{\sqrt{k}}\right)\varrho \big(k^{-\frac{1}{9}}\|\uu\|\big)d\uu dV_{X}\\
=k^{-\frac{d_{P}-1}{2}}\cdot\frac{2^{\frac{d_{G}}{2}}d_{\nu_{G}}^{2}}{(\sqrt{2})^{d_{T}-1}\pi^{d_{T}-1}}\left(\frac{\|\nu_{T}\|k}{\pi}\right)^{d_{M}-\frac{d_{P}-1}{2}}\cdot\\
\int_{X_{0,\nu_{T}}}\int_{\mathbb{R}^{d_{P}-1}}\frac{f(\pi(x))}{\|\Phi_{T}\|^{d_{M}-\frac{d_{P}-1}{2}+1}\mathcal{D}(\pi(x))}e^{-\lambda_{\nu_{T}}2\|\uu\|^2}\varrho \big(k^{-\frac{1}{9}}\|\uu\|\big)d\uu dV_{X}(x) + \cdots,
\end{multlined}
\end{equation}
\noindent
where $d_{P}=d_{G}+d_{T}$ and the dots denote lower order terms.
Now we evaluate the Gaussian integral, that is the same of $(\ref{boiaboiagausss})$ in the Corollary $\ref{cor:quarto}$.
Substituting the result in $(\ref{telufficio57314})$ we obtain the following expression:

\begin{equation*}
\begin{multlined}[t][12.5cm]
\mathfrak{T}\left(T_{\nu_{G},k\nu_{T}}[f]\right)=\frac{d_{\nu_{G}}^{2}}{2^{d_{T}-1}\pi^{d_{T}-1}}\left(\frac{\|\nu_{T}\|k}{\pi}\right)^{d_{M}-d_{P}+1}\cdot\\
\cdot\int_{X_{0,\nu_{T}}}\frac{f(\pi(x))\|\Phi_{T}(\pi(x))\|^{-d_{M}+d_{P}-2}}{\mathcal{D}(\pi(x))}dV_{X}(x)+ \cdots.
\end{multlined}
\end{equation*}

The proof is complete.
\hfill $\Box$

\begin{center}
\textbf{Acknowledgments}
\end{center}

This article is based on my doctoral dissertation at the University of Milano Bicocca under the supervision of Professor Roberto Paoletti. I am grateful to my advisor for introducing me to this area of mathematics and guiding me patiently. I am also endebted to Professor Andrea Loi for carefully reading this manuscript and suggesting many significants improvements in exposition and organization.

\end{document}